\documentclass[12pt]{amsart}

\usepackage{amssymb}
\usepackage{mathrsfs}
\usepackage{graphicx}
\usepackage{psfrag}

\def\makemargins{
	\oddsidemargin .25in
	\evensidemargin .25in
	\textwidth 6.00in
}
\makemargins
\theoremstyle{plain}
\newtheorem{theorem}[subsubsection]{Theorem}
\newtheorem{lemma}[subsubsection]{Lemma}
\newtheorem{proposition}[subsubsection]{Proposition}
\newtheorem{corollary}[subsubsection]{Corollary}

\theoremstyle{definition}
\newtheorem{definition}[subsubsection]{Definition}

\theoremstyle{remark}

\newcommand{\bbC}{{\mathbb C}}
\newcommand{\bbA}{{\mathbb A}}
\newcommand{\bbF}{{\mathbb F}}
\newcommand{\bbG}{{\mathbb G}}
\newcommand{\bbP}{{\mathbb P}}
\newcommand{\bbZ}{{\mathbb Z}}
\newcommand{\bbQ}{{\mathbb Q}}

\newcommand{\cE}{{\mathcal E}}
\newcommand{\cT}{{\mathcal T}}
\newcommand{\cH}{{\mathcal H}}

\newcommand{\sS}{{\mathscr S}}

\newcommand{\sets}[1]{[\![#1]\!]}
\renewcommand{\tilde}{\widetilde}
\newcommand{\card}[1]{|#1|}

\renewcommand{\emptyset}{\varnothing}

\newcommand{\Grass}{{\mathbf G}}
\newcommand{\Flag}{{\mathbf F}}
\newcommand{\SL}{\operatorname{SL}}
\newcommand{\spl}{\text{split}}
\newcommand{\shi}{\text{shift}}

\DeclareMathOperator{\cod}{codim}

\catcode`\@=11
\def\@secnumfont{\bfseries}
\catcode`\@12

\title{Geometry of the tetrahedron space}

\newif \ifdraft

\def \makeauthor{
\author{Eric Babson}
\address{Department of Mathematics\\
University of Washington\\
Seattle, WA  98195}
\email{babson@math.washington.edu}

\author{Paul E. Gunnells}
\address{Department of Mathematics and Statistics\\
University of Massachusetts\\
Amherst, MA  01003}
\email{gunnells@math.umass.edu}

\author{Richard Scott}
\address{Department of Mathematics and Computer Science\\
Santa Clara University\\
Santa Clara, CA  95053}
\email{rscott@math.scu.edu}
}
\makeauthor
\begin{document}

\thanks{Research partially supported by the NSF}
\begin{abstract}
Let $X^{\circ }$ be the space of all labeled tetrahedra in $\bbP^{3}$.
In \cite{BGS} we constructed a smooth symmetric compactification
$\tilde{X}$ of $X^{\circ }$.  In this article we show that the
complement $\tilde{X}\smallsetminus X^{\circ }$ is a divisor with
normal crossings, and we compute the cohomology ring
$H^*(\tilde{X};\bbQ)$.
\end{abstract}

\maketitle

\section{Introduction}
In this article we describe the natural stratification and
intersection theory of the space $\tilde{X}$ of \emph{complete
tetrahedra}, a complex projective variety we constructed in \cite{BGS}
that provides a natural compactification of the variety $X^\circ$ of
nondegenerate tetrahedra in $\bbP^3$.  The importance of $\tilde{X}$
lies in both its connection to work of Schubert, who described an
analogous space for triangles in $\bbP^{2}$ \cite{Schubert} (see also
\cite{Collino-Fulton,RS1,RS2,RS3}), and its relation to more
recently constructed compactifications of configuration varieties
\cite{FM,Keel,Magyar}.  In particular, $\tilde{X}$ provides a natural
setting for studying certain enumerative questions, and for studying
generalized Schur modules for $GL_{4}$ (cf. \cite{vdk-mag}).

There are several desired properties that guide the construction of a
compactification $\tilde{V}$ of a configuration variety $V^{\circ}$:
\begin{itemize}
\item $\tilde{V}$ should be smooth;
\item any group action on $V^{\circ}$ should extend to an action on
$\tilde{V}$;
\item the complement $\tilde{V}\smallsetminus V^{\circ}$ should have a natural
stratification (ideally it should be a divisor with normal crossings); and
\item the cohomology (or intersection) ring of $\tilde{V}$
should have an explicit description in terms of the classes of
closures of strata.
\end{itemize}
The first two properties were verified for the space $\tilde{X}$ of complete
tetrahedra in \cite{BGS}. In the present paper, we show that the
last two properties hold as well.

\begin{figure}[htb]
\psfrag{1}{$\scriptscriptstyle{1}$}
\psfrag{2}{$\scriptscriptstyle{2}$}
\psfrag{3}{$\scriptscriptstyle{3}$}
\psfrag{4}{$\scriptscriptstyle{4}$}
\psfrag{12}{$\scriptscriptstyle{12}$}
\psfrag{13}{$\scriptscriptstyle{13}$}
\psfrag{14}{$\scriptscriptstyle{14}$}
\psfrag{23}{$\scriptscriptstyle{23}$}
\psfrag{24}{$\scriptscriptstyle{24}$}
\psfrag{34}{$\scriptscriptstyle{34}$}
\psfrag{123}{$\scriptscriptstyle{123}$}
\psfrag{124}{$\scriptscriptstyle{124}$}
\psfrag{134}{$\scriptscriptstyle{134}$}
\psfrag{234}{$\scriptscriptstyle{234}$}
\begin{center}
\includegraphics[scale = 0.5]{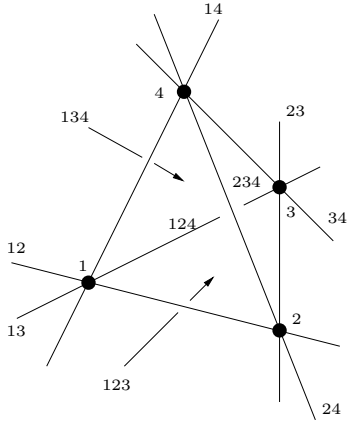}
\end{center}
\caption{A nondegenerate tetrahedron.\label{fig:nondegenerate}}
\end{figure}

To provide some insight into the space $\tilde{X}$, we first describe
the variety of nondegenerate tetrahedra $X^{\circ}$, its canonical
singular compactification $X$, and the divisor at infinity
$X\smallsetminus X^{\circ}$.  Given $4$ general points $P_i$ ($1\leq
i\leq 4$) in $\bbP^3$, each pair determines a line $P_{ij}$, and each
triple determines a plane $P_{ijk}$.  Thus, the $P_{i}$ determine a
point in the product
\begin{equation}\label{eq:product}
(\bbP^3)^4\times\bbG (2,4)^6\times(\hat\bbP^3)^4,
\end{equation}
where $\bbG (2,4)$ is the Grassmannian of lines in $\bbP^{3}$, and
$\hat\bbP^{3}$ is the projective space of hyperplanes in $\bbP^{3}$.
Any such point can be thought of as a nondegenerate tetrahedron, in
the sense that it corresponds to a collection of subspaces in
$\bbP^{3}$ arranged to form a tetrahedron
(Figure~\ref{fig:nondegenerate}).  We let $X^{\circ}$ be the
subvariety of \eqref{eq:product} consisting of all such points.

\begin{figure}[htb]
\begin{center}
\includegraphics[scale = .4]{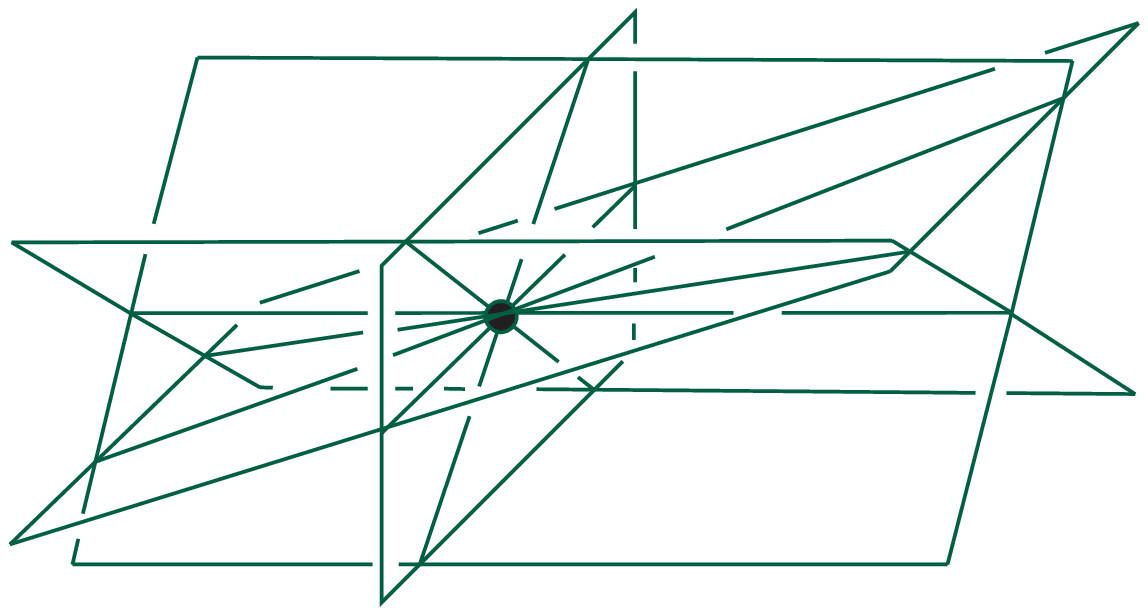}\hspace{.2in}
\raisebox{.1in}{\includegraphics[scale = .4]{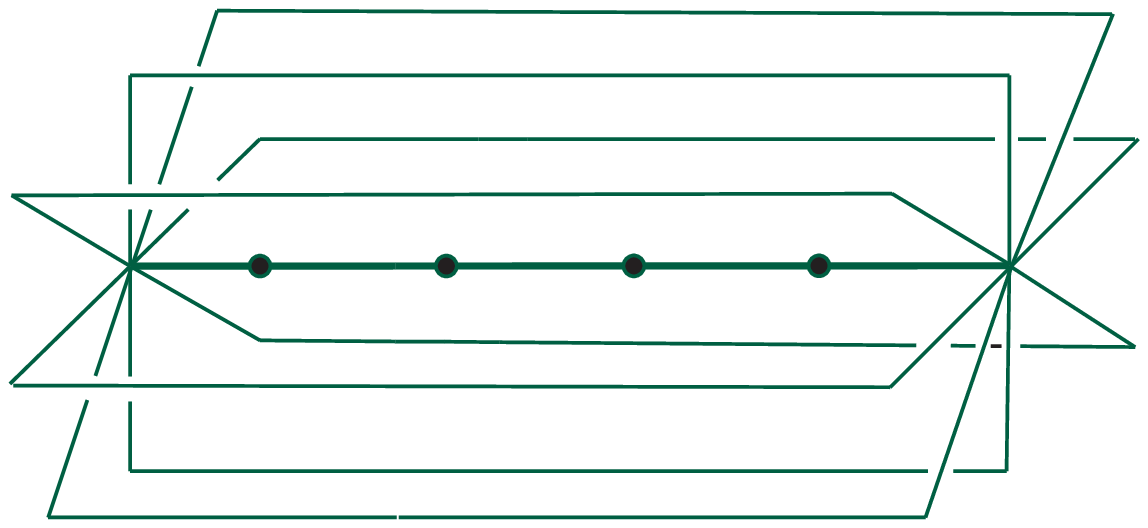}}\hspace{.2in}
\raisebox{.3in}{\includegraphics[scale = .4]{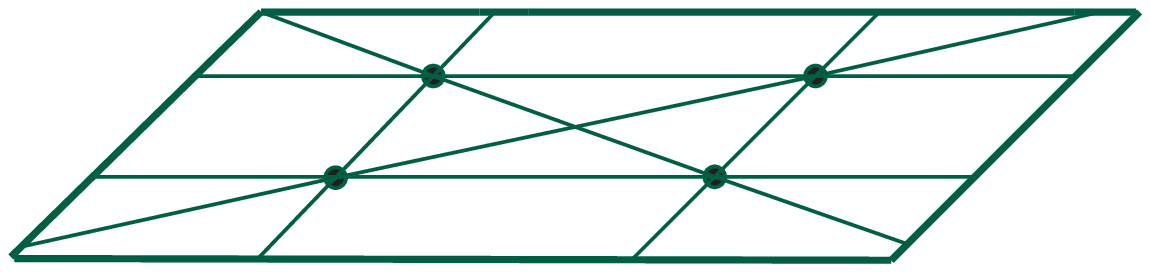}}
\end{center}
\caption{\label{fig:shifting-div} The shifting divisors $A$, $B$, $A^*$.}
\end{figure}

\begin{figure}[htb]
\begin{center}
\includegraphics[scale = .55]{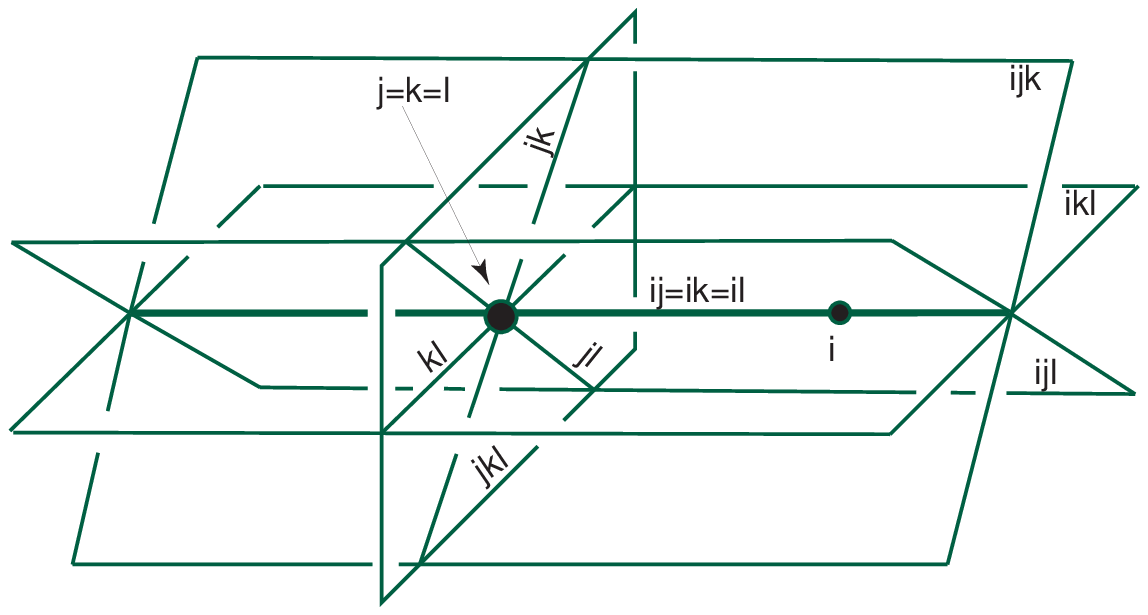}\hspace{.3in}
\raisebox{.3in}{\includegraphics[scale = .55]{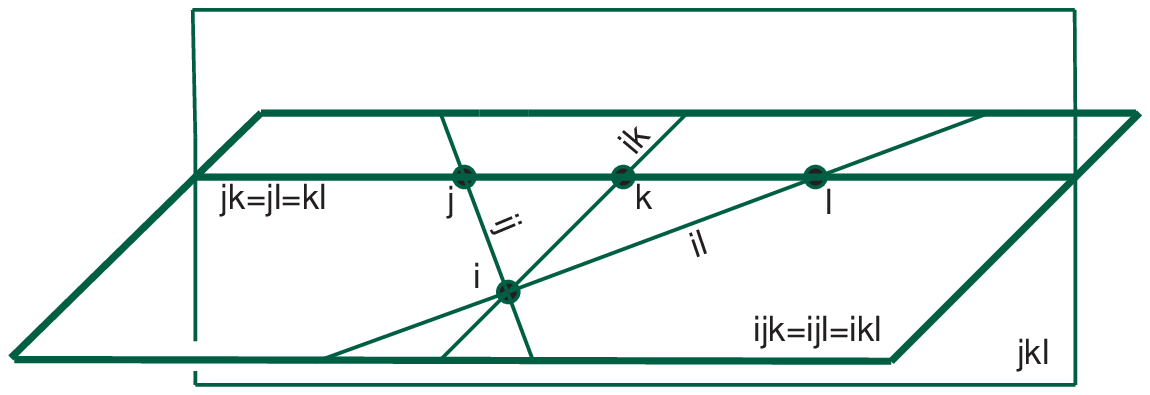}}\\
\includegraphics[scale = .55]{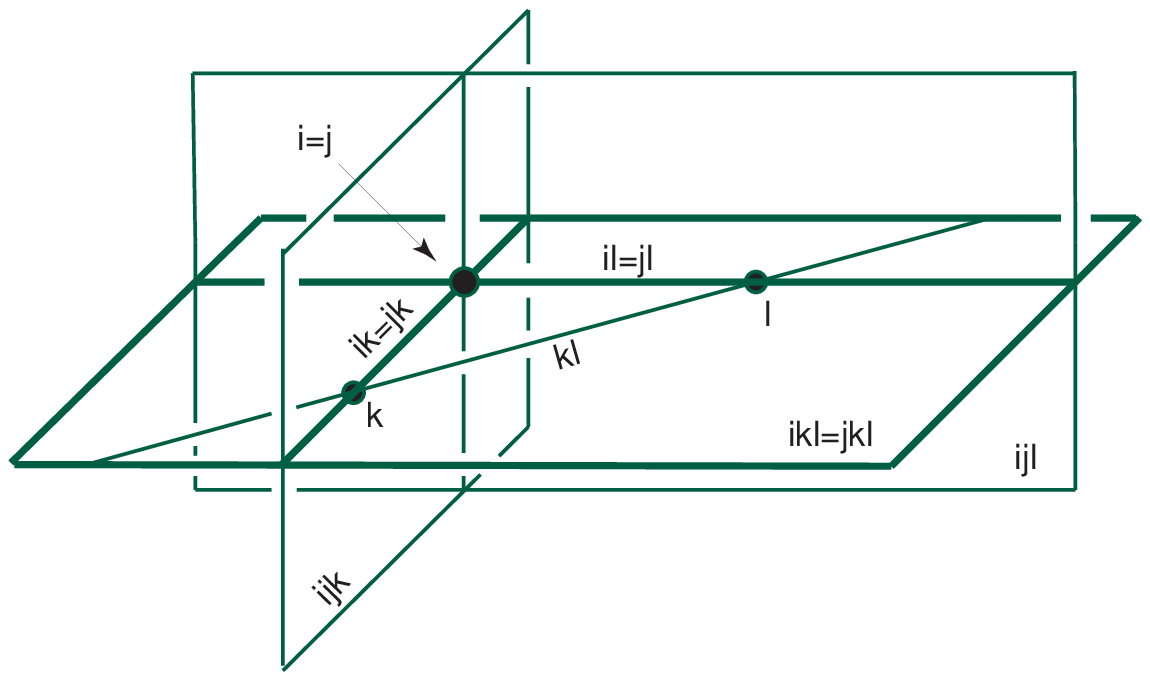}\hspace{.3in}
\raisebox{.3in}{\includegraphics[scale = .55]{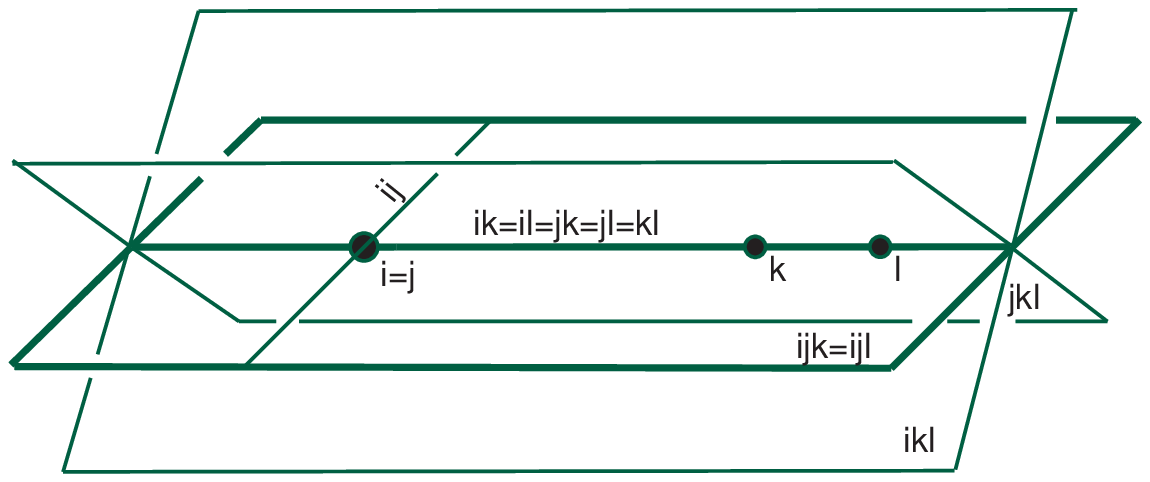}}
\end{center}
\caption{\label{fig:split-div} The split divisors $C_i$, $C_i^*$,
$D_{ij}$, and $E_{ij}$.}
\end{figure}

By taking the closure of $X^{\circ}$ in \eqref{eq:product}, one
obtains the {\em canonical space of tetrahedra} $X$.  The points in
$X\smallsetminus X^{\circ}$ are \emph{degenerate tetrahedra}; that is
to say, they correspond to configurations of subspaces that can be
obtained as limits of families of nondegenerate tetrahedra.  It turns
out that, up to symmetry, there are $7$ maximal combinatorial types of
degenerate tetrahedra parameterized by $X$
(Figures~\ref{fig:shifting-div} and~\ref{fig:split-div}).  These types
are maximal in the sense that any degenerate tetrahedron in
$X\smallsetminus X^{\circ}$ either has one of these types or can be
obtained by a further degeneration of one of these types.  If one
labels these configurations with subsets of $\{1, 2, 3, 4 \}$ as in
Figure \ref{fig:nondegenerate}, one finds altogether $23$ maximal
combinatorial types (namely, $A,B,A^*,C_i,C_i^*,D_{ij},E_{ij}$).  The
closure of the locus of all points in $X$ of a given type is an
irreducible divisor in $X$, and the union of these divisors is
precisely the complement of $X^{\circ}$ in $X$.

The variety $X$ is singular, and the locus of degenerate tetrahedra is
relatively complicated.  In \cite{BGS}, we constructed a smooth
symmetric compactification $\tilde{X}$ of $X^{\circ}$ that dominates
$X$; one obtains $\tilde{X}$ by taking the closure of $X^{\circ}$ in a
larger ambient space than \eqref{eq:product} (see
\ref{ss:bundle-def}).  The variety $\tilde{X}$ has the property that
the fiber of the map $\tilde{X}\rightarrow X$ over a degenerate
tetrahedron of one of the $23$ types described above is a single
point.  Hence, we obtain a collection of irreducible divisors in
$\tilde{X}$ by taking proper transforms of the divisors described
above.  Our first result (Theorem~\ref{thm:normalcrossings}) shows
that in the resolution $\tilde{X}\rightarrow X$, no new divisors are
introduced, and that furthermore one of our criteria for a good
compactification of $X^{\circ}$ is met:

\smallskip
\noindent{\bf Theorem.}
{\em Let $Z\subset\tilde{X}$ be the union of the $23$ irreducible
divisors.  Then $Z=\tilde{X}\smallsetminus X^{\circ}$, and is a
divisor with normal crossings.}
\smallskip

For each of the $4$ points, $6$ lines, and $4$ planes in a tetrahedron, there
is a projection from $\tilde{X}$ to the corresponding
Grassmannian; this is simply the composition of $\tilde{X}\rightarrow
X$ with projection to the appropriate factor of \eqref{eq:product}.
By pulling back the special Schubert varieties \cite[p. 271]{fulton} via
these maps, we obtain additional divisors $Y_{I}$ in $\tilde{X}$,
where $I$ ranges over all proper nonempty subsets of $\{1,2,3,4
\}$.
For each $I$, the divisor $Y_{I}$ consists of those points in
$\tilde{X}$ whose image tetrahedra have their subspace labeled $I$
meeting a given codimension-$\card{I}$ subspace $W\subset \bbP^{3}$.  We call
these divisors \emph{special position divisors}.  Together with the
$23$ divisors mentioned above, they generate the cohomology ring of
$\tilde{X}$.  In fact, the complete ring structure is given by the
following presentation (Theorem~\ref{thm:ringpres}):

\smallskip
\noindent{\bf Theorem.}
{\em The cohomology ring $H^*(\tilde{X};\bbQ)$ is generated in degree
$2$ by the Poincar\'e duals of the $23$ divisors from
Figures~\ref{fig:shifting-div} and~\ref{fig:split-div}, and the
special position divisors $Y_{I}$.  If we denote these
dual classes by $a$, $b$, $a^*$, $c_i$, $c_i^*$, $d_{ij}$, $e_{ij}$,
$y_i$, $y_{ij}$, $y_{ijk}$ (where the indices denote {\em unordered}
subsets of $\{1,2,3,4\}$), then the ideal of relations is generated by
the following polynomials:
\begin{enumerate}
\item[(i)] $y_{ij}-y_i-y_j+a+c_k+c_l+d_{ij}+e_{ij}$,\\
$y_{ijk}-y_{ij}-y_{ik}+y_i+b+c_i+c_l^*+d_{jk}+e_{jk}+e_{il}+e_{jl}+e_{kl}$,\\
$y_{ij}-y_{ijk}-y_{ijl}+a^*+c_i^*+c_j^*+d_{kl}+e_{ij}$.\\
\item[(ii)]  $c_ic_j,\;\; c_i^*c_j^*,\;\; c_id_{ij},\;\;c_i^*d_{ij},\;\;
c_ie_{ij},\;\; c_i^*e_{jk},\;\; d_{ij}e_{ik},\;\; e_{ij}e_{ik},\;\;
e_{ij}e_{kl}$.\\
\item[(iii)] $a(y_i-y_j),\; b(y_{ij}-y_{ik}),\; a^*(y_{ijk}-y_{ijl}),\\
c_i(y_j-y_k),\; c_i(y_{ij}-y_{ik}),\; c_i^*(y_{jk}-y_{jl}),\;
c_i^*(y_{ijk}-y_{ijl}),\\
d_{ij}(y_i-y_j),\; d_{ij}(y_{ik}-y_{jk}),\; d_{ij}(y_{ikl}-y_{jkl}),\\
e_{ij}(y_i-y_j),\; e_{ij}(y_{kl}-y_{ik}),\;
e_{ij}(y_{ijk}-y_{ijl})$.\\
\item[(iv)] $y_i^2+y_{ij}^2+y_{ijk}^2-y_iy_{ij}-y_{ij}y_{ijk}$,\\
$y_{ij}^3-2y_iy_{ij}^2+2y_i^2y_{ij}$,\\
$y_i^4$.
\end{enumerate}
}
\smallskip

The relations in this presentation all come from simple geometric
considerations.  Those in (ii), for example, arise because certain
pairs of divisors are disjoint in $\tilde{X}$.  The relations in (i)
come from rational equivalences induced by rational functions
on $\tilde{X}$ corresponding to certain cross-ratios.  More details appear in
\ref{ss:relations}.

The structure of the paper is as follows.  In
Section~\ref{s:background}, we recall the necessary constructions and
results from \cite{BGS}, and state other results we need whose proofs
are easy adaptations of results from \cite{BGS}.  The main tool we use
throughout this paper is the local description of $\tilde{X}$ given in
Theorem~\ref{thm:local-equations}.  In Section
\ref{s:stratifications}, we describe the combinatorial diagrams we use
to decompose $\tilde{X}$ into strata, show that the union of the
strata of codimension $\geq 1$ is a divisor with normal crossings, and
show that this decomposition into subvarieties is indeed a
stratification.  In Section \ref{s:betti} we compute the topological
Betti numbers of $\tilde{X}$ by computing the Hasse-Weil zeta function
of $\tilde{X}$ and using the Weil conjectures.  Finally, in
Section~\ref{s:ring} we compute the cohomology ring $H^{*}
(\tilde{X};\bbQ )$.  A key step in this computation is the use of the
Betti numbers from Section \ref{s:betti} to verify that the list of
relations in the above presentation suffices.

%
\section{Background}\label{s:background}
We recall the setup and basic constructions from \cite{BGS}.  Proofs
of all the statements in this section can be found in \cite{BGS} or
are straightforward generalizations of results in \cite{BGS}.

\subsection{The singular space of tetrahedra}
Let $\sets{4}$ be the set $\{1,2,3,4\}$.  For any proper nontrivial
subset $I\subset\sets{4}$, we let $\Grass_I$ denote the Grassmannian
of $\card{I}$-planes in $\bbC^4$.  Let $\Grass$ be the product
\[\Grass := \prod _{I} \Grass _{I},\]
and for any $x\in\Grass$, let $x_I\in\Grass_I$ be its image in
the $I$th factor.

Let $e_1,\ldots,e_4$ be the standard basis for $\bbC^4$.  For any
proper nontrivial subset $I\subset\sets{4}$, let $e_I\in\Grass_I$
be the subspace spanned by $\{e_i\;|\;i\in I\}$, and let
$e\in\Grass$ be the point $e=(e_I)_{I\subset\sets{4}}$.  The group
$G=\SL_{4}(\bbC)$ acts diagonally on $\Grass$.

\begin{definition}
The {\em space of nondegenerate tetrahedra}, which we denote by
$X^{\circ}$, is the orbit $G\cdot e\subset\Grass$.  The {\em
(canonical) space of tetrahedra}, denoted by $X$, is its closure
$\overline{G\cdot e}\subset\Grass$.  A point in $X^{\circ}$ is a {\em
nondegenerate tetrahedron}, a point in $X$ is a {\em tetrahedron},
and a point in $X\smallsetminus X^{\circ}$ is a {\em degenerate
tetrahedron}.
\end{definition}

Let $\Flag$ be the variety of full flags in $\bbC^{4}$.

\begin{proposition}\hspace{1in}
\begin{itemize}
\item The $G$-action on $\Grass$ restricts to an action on $X$.
\item The symmetric group $S_4$ acts on $X$ via the action on $\Grass$
induced from the natural $S_{4}$ action on $\sets{4}$.
\item Each projection $X\rightarrow\Grass_I$ is a $G$-equivariant
locally trivial fibration.
\item The projection
$X\rightarrow\Flag\subseteq\Grass_1\times\Grass_{12}\times\Grass_{123}$ is a
$G$-equivariant locally trivial fibration.
\end{itemize}
\end{proposition}

\subsection{The smooth space of tetrahedra}\label{ss:bundle-def}
\begin{figure}[htb]
\begin{center}
\psfrag{1}{\hspace{-.03in}\scriptsize$1$}
\psfrag{2}{\hspace{-.05in}\scriptsize$2$}
\psfrag{3}{\scriptsize$3$}
\psfrag{4}{\hspace{-.03in}\scriptsize$4$}
\psfrag{12}{\hspace{-.05in}\scriptsize$12$}
\psfrag{13}{\scriptsize$13$}
\psfrag{14}{\hspace{-.05in}\scriptsize$14$}
\psfrag{23}{\hspace{-.05in}\scriptsize$23$}
\psfrag{24}{\hspace{-.05in}\scriptsize$24$}
\psfrag{34}{\scriptsize$34$}
\psfrag{123}{\hspace{-.05in}\scriptsize$123$}
\psfrag{124}{\scriptsize$124$}
\psfrag{134}{\hspace{-.05in}\scriptsize$134$}
\psfrag{234}{\hspace{-.05in}\scriptsize$234$}
\includegraphics[scale = .7]{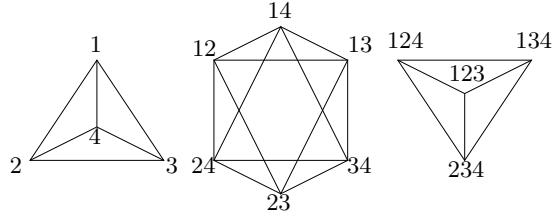}
\end{center}
\caption{The three $3$-dimensional hypersimplices.\label{fig:diagram}}
\end{figure}

Let $\Delta _{1}$, $\Delta_{2} $, $\Delta _{3}$ be the $3$-dimensional
hypersimplices, with vertices indexed by proper nontrivial subsets
$I\subset\sets{4}$ (Figure~\ref{fig:diagram}).  For each such $I$, we
let $F_I\rightarrow X$ be the pull-back of the tautological bundle
on $\Grass_I$.

The edges of the hypersimplex $\Delta_i$ are indexed by pairs
$\{I,J\}$ of $i$-element subsets of $\sets{4}$ satisfying $\card{I\cap J}=i-1$ and
$\card{I\cup J}=i+1$.  Let $\cE$ be the set of all edges of the
hypersimplices.  For each $\alpha=\{I,J\}\in\cE$, we let
$P_{\alpha}\rightarrow X$ be the $\bbP^1$-bundle defined by
$\bbP(F_{I\cup J}/F_{I\cap J})$.   Each $P_{\alpha}$ has canonical
sections $s_I$ and $s_J$ defined by $s_I(x)=x_I/x_{I\cap J}$ and
$s_J(x)=x_J/x_{I\cap J}$.   We let $E_{\alpha}\rightarrow X$ denote the
$\bbP^1\times\bbP^1$-bundle $P_{\alpha}\times_X P_{\alpha}$ with
diagonal subbundle $D_{\alpha}\rightarrow X$.  The sections $s_I$ and
$s_J$ determine a section $s_I\times s_J$ of $E_{\alpha}$ which we
denote by $s_{\alpha}$.  The significance of the
section $s_{\alpha}$ is that the preimage of the diagonal $D_{\alpha}$
is precisely the locus where the $I$th and $J$th face of a tetrahedron
coincide, i.e., $s_{\alpha}(x)\in D_{\alpha}$ if and only if
$x_I=x_J$.

\begin{figure}[htb]
\begin{center}
\includegraphics[scale = .5]{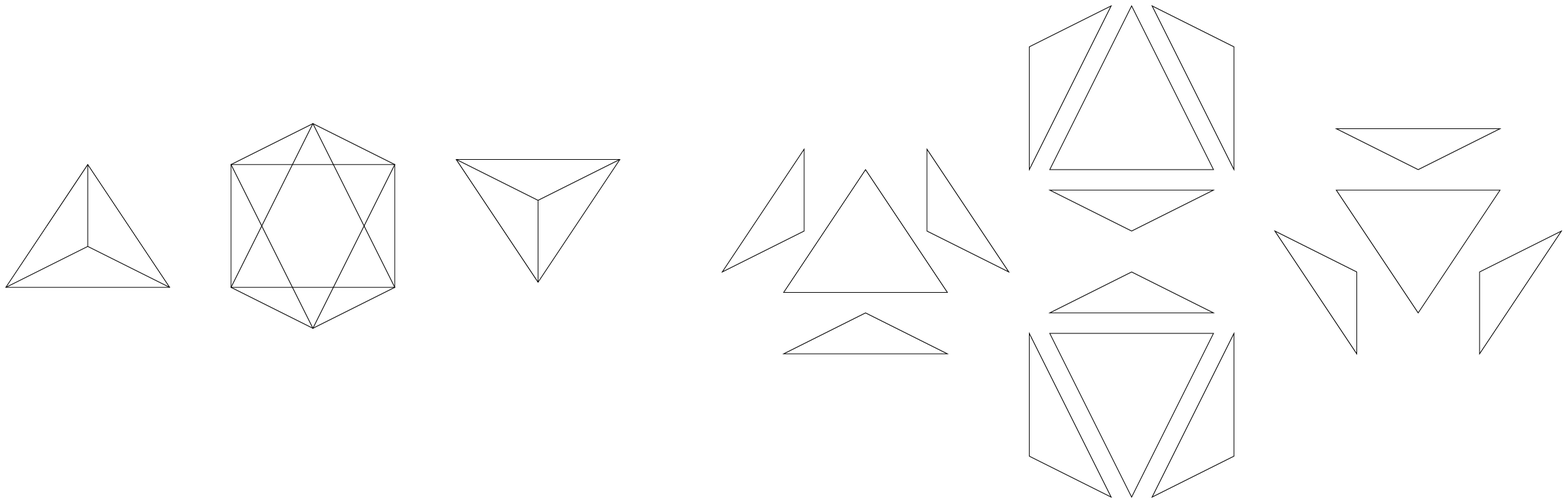}
\end{center}
\caption{The set $\cH$.\label{fig:fiber}}
\end{figure}

Let $\cH$ be the set of all hypersimplex faces of dimension $\geq 2$.
This set consists of the three $3$-dimensional hypersimplices and
their $16$ triangular faces (Figure~\ref{fig:fiber}).  For each
$\beta\in\cH$, we define $E_{\beta}\rightarrow X$ to be the product
bundle $\prod E_{\alpha}$, where $\alpha$ ranges over edges of $\beta$.
We let $D_{\beta}$ be the corresponding product $\prod D_{\alpha}$ of
diagonals, and we let $s_{\beta}$ be the section $\prod s_{\alpha}$ of
$E_{\beta}\rightarrow X$.  As above, we have $s_{\beta}(x)\in
D_{\beta}$ if and only if $x_I=x_J$ for all vertices $I,J\in\beta$.

For each $\beta\in\cH$, let $(E_{\beta})_{\#}\rightarrow X$ be
the blow-up of $E_{\beta}$ along the subbundle $D_{\beta}$, and let
$E_{\#}\rightarrow X$ be the product
\[E_{\#}=\prod_{\beta\in\cH}(E_{\beta})_{\#}.\]
The product of the sections $s_{\beta}$ induces a rational section
$s:X\dashrightarrow E_{\#}$ whose restriction to the open set
$X^{\circ}$ is regular.  We can now define the variety of interest in
this paper.

\begin{definition}
The {\em space of complete tetrahedra}, which we denote by
$\tilde{X}$, is the closure of the image of the rational section $s$, i.e.,
$\tilde{X}=\overline{s(X^{\circ})}$.  A point $\tilde{x}\in\tilde{X}$
is a {\em complete tetrahedron}.
\end{definition}

\begin{proposition}\hspace{1in}
\begin{itemize}
\item $\tilde{X}$ is a smooth projective variety.
\item There is a natural surjective birational morphism
$\tilde{X}\rightarrow X$ defined by restricting the bundle projection
$E_{\#}\rightarrow X$.
\item There are natural actions of $G$ and the symmetric group $S_4$
on $\tilde{X}$, and the projection $\tilde{X}\rightarrow X$ is
equivariant with respect to both.
\item The compositions $\tilde{X}\rightarrow X\rightarrow\Grass_I$
and $\tilde{X}\rightarrow X\rightarrow\Flag$ are locally trivial
$G$-equivariant fibrations.
\end{itemize}
\end{proposition}

\subsection{Local equations for $X$ and $\tilde{X}$}\label{ss:local-equations}

For any flag $V:=(V_1\subset V_2\subset V_3)\in\Flag$, let
$\Flag(V)\subset\Flag$ denote the open set consisting of flags in
general position to $V$.  Thus, $\Flag(V)$ is a $6$-dimensional affine
space.  For any $x\in X$ (respectively, $\tilde{x}\in\tilde{X}$), we
call its image $x_I\in\Grass_I$ the {\em $I$th plane of $x$ (resp.,
$\tilde{x}$)}.  Let $U(V)$ (resp., $=\tilde{U}(V)$) be the open subset
of $X$ (resp., $\tilde{X}$) consisting of those $x$ (resp.,
$\tilde{x}$) such that $x_I$ is in general position to $V$ for each
$I\subset\sets{4}$.  By varying $V$ we obtain open covers
$\{U(V)\}_{V\in\Flag}$ and $\{\tilde{U}(V)\}_{V\in\Flag}$ for $X$ and
$\tilde{X}$, respectively \cite[Lemma
4.3]{BGS}.  In \cite[5.11]{BGS}, we defined closed embeddings
\[U(V)\subset\Flag(V)\times\bbC^{24}\]
and
\[\tilde{U}(V)\subset\Flag(V)\times\bbC^{24}\times(\bbP^5)^2\times
\bbP^{11}\times(\bbP^2)^{16},\]
and determined defining equations for their images.  We recall these
equations in this subsection.

Let $\bbA_{\cE}$ be the affine space with coordinates
$\{u_{\alpha}\;|\;\alpha\in\cE\}$, and for each $\beta\in\cH$, let
$\bbP_{\beta}$ denote the projective space with coordinates
$\{u_{\alpha,\beta}\;|\;\alpha\in\cE,\;\alpha\subset\beta\}$.
We let $\cE_{\#}$ denote the set of all pairs $(\alpha,\beta)$ such
that $\alpha\in\cE$ and $\alpha\subset\beta$.  Since the elements of
$\cE$ appear as the edges in Figure~\ref{fig:diagram} and the elements
of $\cE_{\#}$ appear as the edges in Figure~\ref{fig:fiber}, we shall
often refer to any elements of $\cE$ or $\cE_{\#}$ as an {\em edge}.
For any edge $\alpha$ or $(\alpha,\beta)$, we call $\alpha\in\cE$ its
{\em original edge}.  We note that each edge corresponds to exactly
one of the coordinates defined above, with the original edges
corresponding to the affine coordinates.  To describe the necessary
equations among these coordinates, we appeal to the geometry of the
specific planar representations of the hypersimplices and their faces
appearing in Figures~\ref{fig:diagram} and \ref{fig:fiber}.

\begin{figure}[htb]
\begin{center}
\includegraphics[scale = .47]{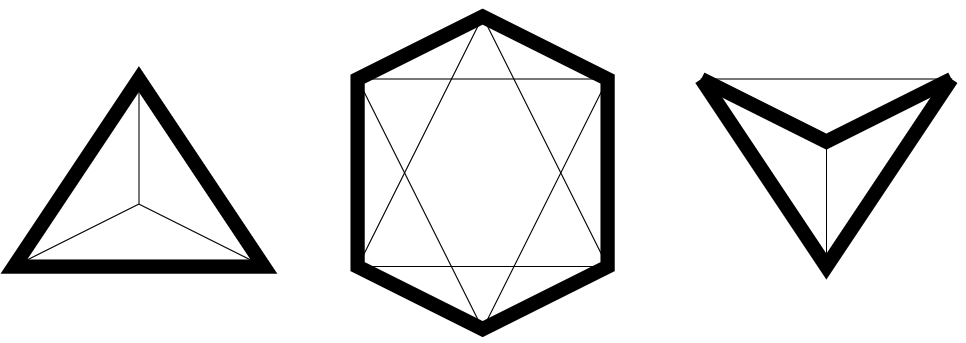}\hspace{.4in}
\includegraphics[scale = .37]{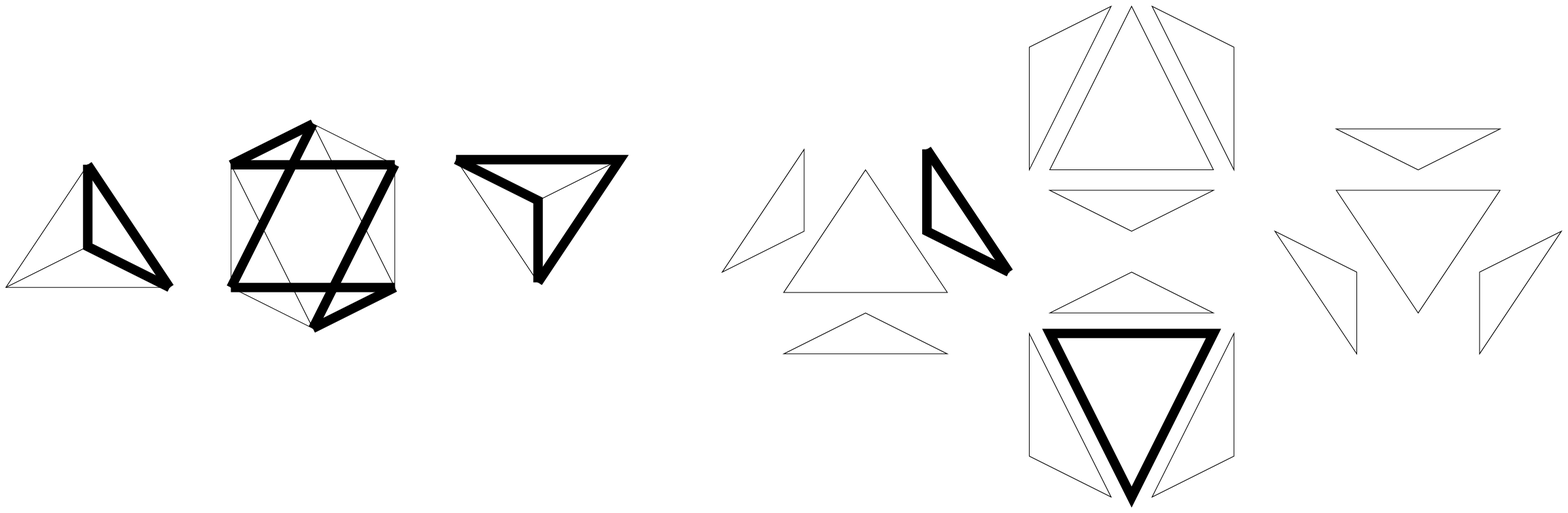}
\end{center}
\caption{\label{fig:polygons} Examples of triangles, quadrilaterals,
and hexagons.}
\end{figure}

\begin{definition}
Let $C=(c_1,\ldots,c_k)$ be a sequence of edges that forms a circuit
in Figure~\ref{fig:diagram} or Figure~\ref{fig:fiber}.  We shall call
$C$ a {\em triangle} if it has length $3$; a {\em
quadrilateral} if it has length $4$ and its original edges
are in $\Delta_1$ or $\Delta_3$; and a {\em hexagon} if it has length
$6$, its original edges are in $\Delta_2$, and its opposite edges are
parallel.  (See Figure~\ref{fig:polygons}.)
\end{definition}

\begin{figure}[htb]
\begin{center}
\includegraphics[scale = .45]{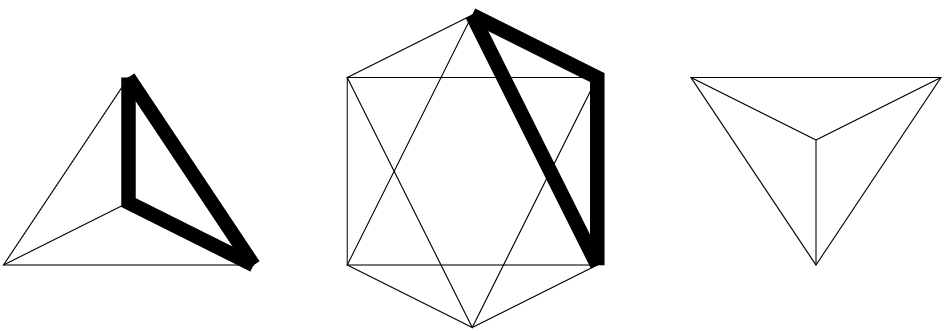}\hspace{.4in}
\includegraphics[scale = .35]{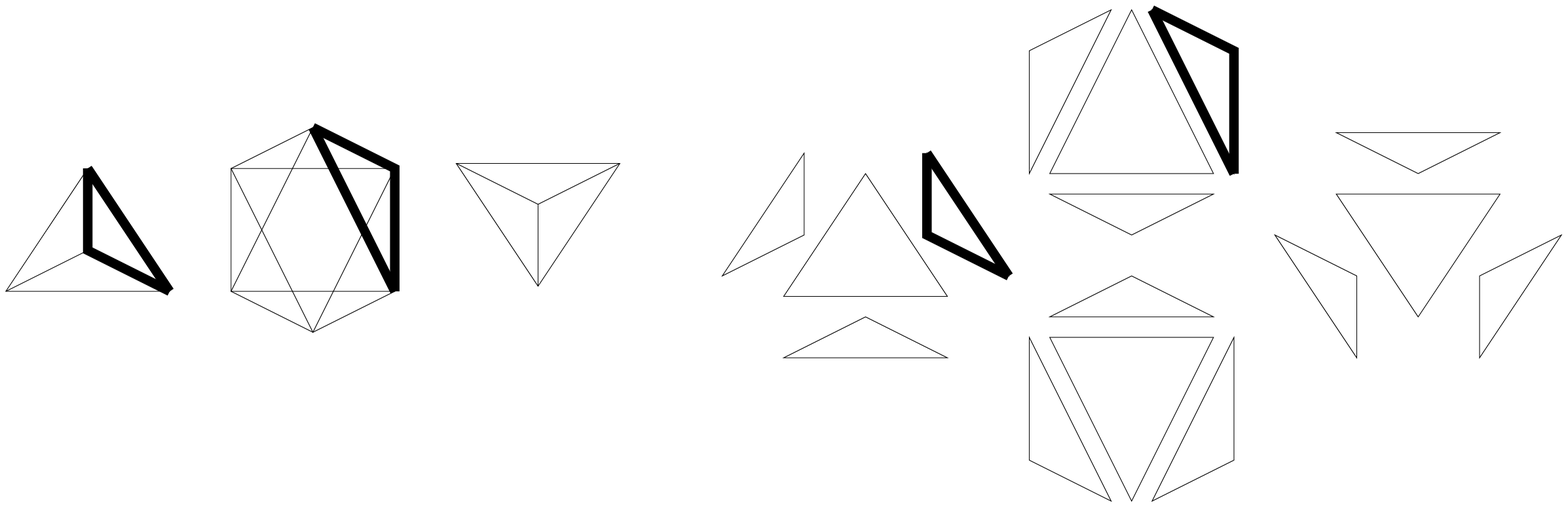}
\end{center}
\caption{\label{fig:rel-tri} Six triangles determining $\binom{6}{2}=15$ pairs
of related triangles.}
\end{figure}

\begin{figure}[htb]
\begin{center}
\includegraphics[scale = .45]{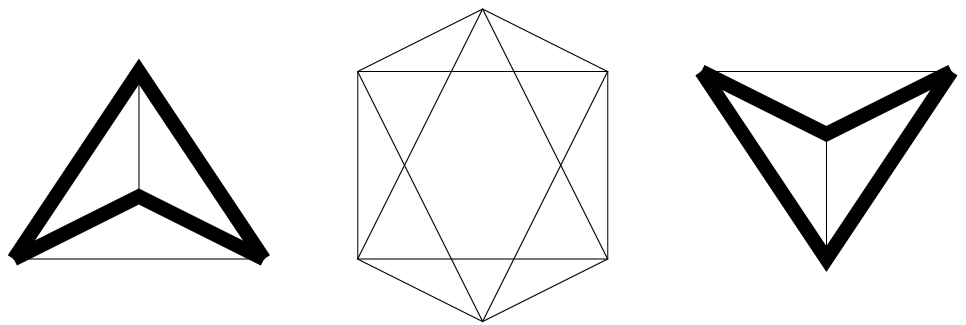}\hspace{.4in}
\includegraphics[scale = .35]{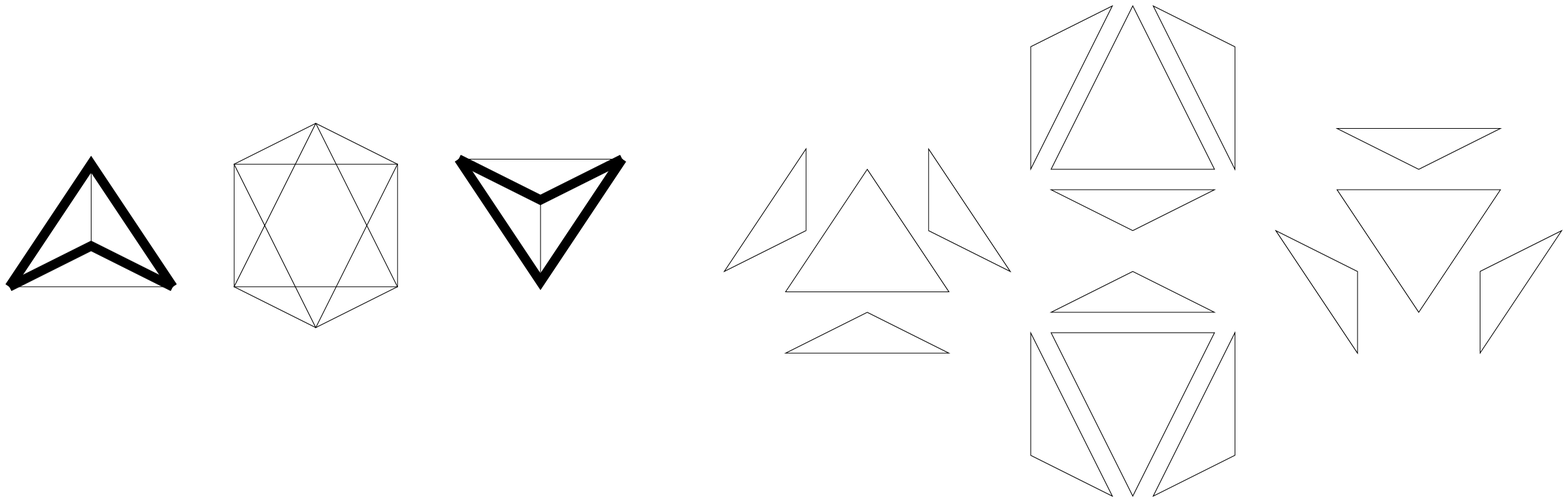}
\end{center}
\caption{\label{fig:rel-quad} Four quadrilaterals determining $4$
pairs of related quadrilaterals.}
\end{figure}

\begin{definition}
Two triangles are {\em related} if either (1) they have the same set
of original edges or (2) their original edges occur in consecutive
$\Delta_i$'s and their images in Figure~\ref{fig:diagram} differ by a
$180^{\circ}$ rotation.  (See Figure~\ref{fig:rel-tri}.)
Corresponding pairs of edges in related triangles will be called {\em
related angles}.  Two quadrilaterals are {\em related} if their
original edges are in $\Delta_1$ and $\Delta_3$ and their images
differ by a $180^{\circ}$ rotation.  (See Figure~\ref{fig:rel-quad}.)
\end{definition}

\begin{theorem}\label{thm:local-equations}
(\cite[4.10 and 5.11]{BGS})
For any flag $V\in\Flag$, let $U=U(V)$ and $\tilde{U}=\tilde{U}(V)$ be
as above.  Then there exist closed embeddings
\[
U\subset\Flag(V)\times\bbA_{\cE}\hspace{.5in}\mbox{and}\hspace{.5in}
\tilde{U}\subset\Flag(V)\times\bbA_{\cE}\times\prod_{\beta\in\cH}
\bbP_{\beta}
\]
whose images are defined set-theoretically by the following equations:
\begin{enumerate}
\item[(1)] $u_a-u_b+u_c=0$ for any triangle $(a,b,c)$.\footnote{There
is a sign convention for our coordinates that depends on an ordering
of the subsets of $\sets{4}$.  We refer to \cite{BGS} for details.}
\item[(2)] $u_au_{a'}=u_bu_{b'}$ for any related angles $(a,b)$
and $(a',b')$.
\item[(3)] $u_au_bu_c=u_{a'}u_{b'}u_{c'}$ for any hexagon
$(a,a',b,b',c,c')$.
\item[(4)] $u_au_{b'}u_cu_{d'}=u_{a'}u_{b}u_{c'}u_{d}$
for any related quadrilaterals $(a,b,c,d)$ and $(a',b',c',d')$.
\end{enumerate}
Moreover, with respect to these embedding, the fibrations
$X\rightarrow\Flag$ and $\tilde{X}\rightarrow\Flag$ are given by
projection to $\Flag(V)$, and the map $\tilde{X}\rightarrow X$ is
given by projection to $\Flag (V)\times \bbA_{\cE}$.
\end{theorem}

For later computations, it will be convenient to write the ambient
variety for $\tilde{U}$ in Theorem~\ref{thm:local-equations} as
\[
\Flag(V)\times\bbA_{\cE_{1}}\times \bbA_{\cE_{2}}\times
\bbA_{\cE_{3}}\times\bbP_{\Delta_{1}}\times \bbP_{\Delta_{2}}\times
\bbP_{\Delta_{3}}\times \prod_{\beta\in\cT} \bbP_{\beta},
\]
where $\cE_{i}$ is the set of edges in $\Delta_{i}$, and where $\cT$
is the set of triangular faces of the $\Delta_{i}$.  Thus, each factor
(except $\Flag(V)$) corresponds to a connected component in
Figure~\ref{fig:polygons}.

\section{Stratifications and Normal crossings}\label{s:stratifications}
Recall that a stratification of a variety $X$ is a collection
$\{X_S\}_{S\in\sS}$ of locally closed subvarieties indexed by a poset
$\sS$ such that
\begin{itemize}
\item $X=\coprod X_{S}$ (disjoint union), and
\item $\overline{X}_S=\coprod_{T\leq S} X_T$.
\end{itemize}
In this section we describe a natural stratification of $\tilde{X}$,
and prove that the closures of the codimension one strata form a
divisor with normal crossings.

\subsection{Diagrams}
We begin by describing the combinatorial data we use to index the strata.
Let $x$ be an element of the canonical space of tetrahedra $X$.
Recall from \ref{ss:bundle-def} that for any edge
$\alpha=\{I,J\}\in\cE$, the image $s_{\alpha}(x)$ will be in the
diagonal $D_{\alpha}$ if and only if the planes $x_I$ and $x_J$
coincide.  This leads us to consider the subset $S(x)\subset\cE$
defined by
\[S(x)=\{\alpha\;|\; s_{\alpha}(x)\in D_{\alpha}\}.\]

\begin{figure}[htb]
\begin{center}
{
\psfrag{1}{\tiny$1$}
\psfrag{2}{\hspace{-.05in}\tiny$2$}
\psfrag{3}{\tiny$3$}
\psfrag{4}{\hspace{-.02in}\vspace*{.05in}\small$_4$}
\psfrag{12}{\hspace{-.08in}\tiny$12$}
\psfrag{13}{\tiny$13$}
\psfrag{14}{\hspace{-.05in}\tiny$14$}
\psfrag{23}{\tiny$23$}
\psfrag{24}{\hspace{-.05in}\tiny$24$}
\psfrag{34}{\tiny$34$}
\psfrag{123}{\hspace{-.05in}\tiny$123$}
\psfrag{124}{\tiny$124$}
\psfrag{134}{\hspace{-.05in}\tiny$134$}
\psfrag{234}{\hspace{-.05in}\tiny$234$}
\includegraphics[scale = .7]{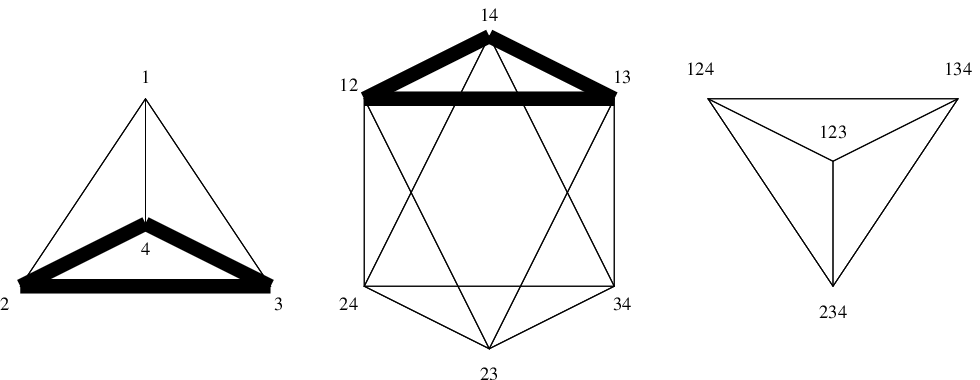}}\hspace{.4in}
\includegraphics[scale = .55]{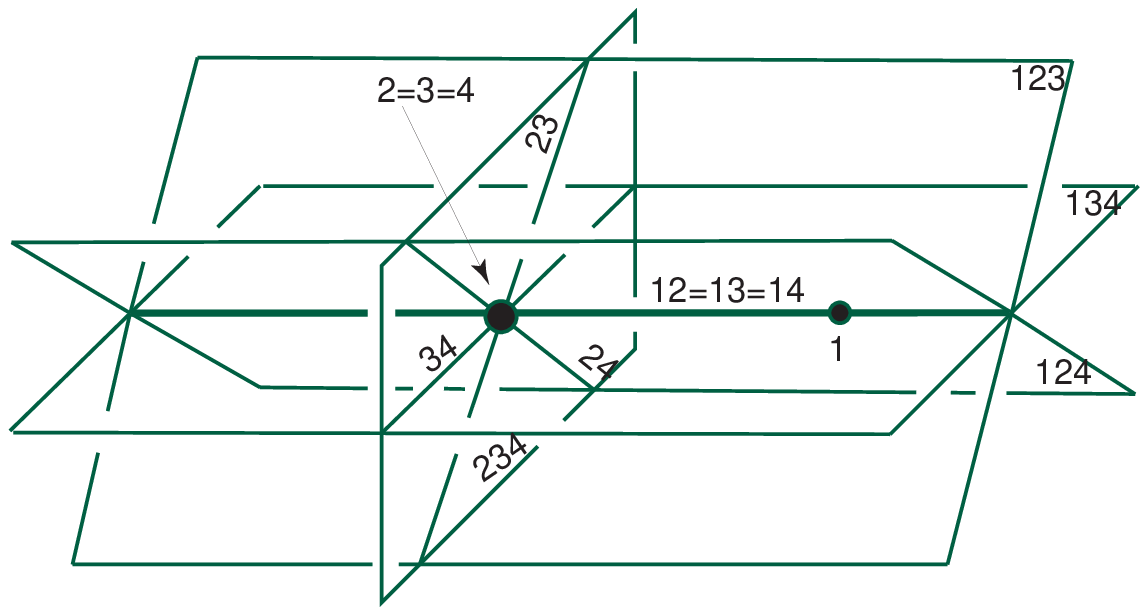}
\end{center}
\caption{\label{fig:c}}
\end{figure}

We represent the subset $S (x)$ graphically by marking in bold the edges in
Figure~\ref{fig:diagram} corresponding to its elements.  These bold
edges encode exactly the projections of $x$ to factors of $\Grass$
that coincide.
For example, if $x$ is a degenerate tetrahedron whose
points $x_2,x_3,x_4$ all coincide, and whose lines
$x_{12},x_{13},x_{14}$ all coincide, and with no other collapsing
among the $x_I$'s, then we have
\[S(x)=\{\{2,3\},\{2,4\},\{3,4\},\{12,13\},\{12,14\},\{13,14\}\}\]
(see Figure~\ref{fig:c}).

Next we shall describe similar marked diagrams for elements of
$\tilde{X}$.  First we give the formal description, which requires
some additional notation.  For any edge $(\alpha,\beta)\in\cE_{\#}$,
we let $E_{\alpha,\beta}\subset E_{\beta}$ be the subbundle
\[E_{\alpha,\beta}=D_{\alpha}\times\prod_{\alpha'\neq\alpha,\; \alpha'\subset\beta}E_{\alpha'}.\]
We let $(E_{\alpha,\beta})_{\#}\rightarrow E_{\alpha,\beta}$ be
its blow-up along
\[D_{\beta}=D_{\alpha}\times\prod_{\alpha'\neq\alpha,\; \alpha'\subset\beta}D_{\alpha'},\]
and we let $D_{\alpha,\beta}$ be the exceptional divisor of this
blow-up.  By functoriality of blow-ups, $(E_{\alpha,\beta})_{\#}$, and
hence $D_{\alpha,\beta}$, are both subvarieties of $(E_{\beta})_{\#}$.

Now suppose $\tilde{x}$ is a point in $\tilde{X}$.  For any
$\beta\in\cH$, we let $\tilde{x}_{\beta}$ denote the image of
$\tilde{x}$ under the projection $E_{\#}\rightarrow(E_{\beta})_{\#}$,
and we define $S_{\#}(\tilde{x})\subset\cE_{\#}$ to be the subset
\[S_{\#}(\tilde{x})=\{(\alpha,\beta)\;|\;\tilde{x}_{\beta}\in
D_{\alpha,\beta}\}.\]

\begin{definition}
A {\em diagram} $\Gamma$ is a pair $(S,S_{\#})$ where $S\subset\cE$
and $S_{\#}\subset\cE_{\#}$.  If $\tilde{x}\in\tilde{X}$ and $x$
is its image in $X$, then the {\em diagram for $\tilde{x}$}, which we
denote by $\Gamma(\tilde{x})$, is the diagram $(S(x),S_{\#}(\tilde{x}))$.
\end{definition}

We represent a diagram graphically by marking in bold the edges in
Figure~\ref{fig:diagram} corresponding to elements of
$S$, and the edges in Figure~\ref{fig:fiber} corresponding to
elements of $S_{\#}$.

To understand the information that a diagram encodes for a point
$\tilde{x}$, consider a curve $\tilde{x}(t)$ in $\tilde{X}$ with
$\tilde{x}(0)=\tilde{x}$ and $x(t)\in X^{\circ}$ for $t\neq 0$
(equivalently, $\tilde{x}(t)=s(x(t))$ for $t\neq 0$).  Then an edge
$\alpha=\{I,J\}$ will be in $S(x)$ if and only if the $I$th and $J$th
planes of $x(t)$ approach each other as $t\rightarrow 0$, i.e.,
\[\lim_{t\rightarrow 0} x(t)_I=\lim_{t\rightarrow 0} x(t)_J.\]
An edge $(\alpha,\beta)$ with $\alpha=\{I,J\}$ will be in
$S_{\#}(\tilde{x})$ if and only if the $I$th and $J$th planes of
$x(t)$ come together and come together {\em faster} than any other
pair $\alpha'=\{I',J'\}$ where $\alpha'$ is an edge of $\beta$ and
$(\alpha',\beta)\not\in S_{\#}(\tilde{x})$.

The following proposition implies that diagrams are compatible with the
local embeddings of Theorem~\ref{thm:local-equations}.  The proof
follows from the discussion in \cite[Section 5]{BGS}.

\begin{proposition}\label{prop:local-diagram}
Let $\tilde{U}\subset\tilde{X}$ be one of the open subvarieties
described in \ref{ss:local-equations}, and let
$\tilde{x}\in\tilde{U}$.  Then the diagram
$\Gamma(\tilde{x})=(S,S_{\#})$ is defined in terms of coordinates by
\[S=\{\alpha\in\cE\;|\;
u_{\alpha}=0\}\hspace{.5in}\mbox{and}\hspace{.5in}
S_{\#}=\{(\alpha,\beta)\in\cE_{\#}\;|\; u_{\alpha,\beta}=0\}.\]
\end{proposition}

\subsection{Classification of diagrams}
For most subsets $S\subset\cE$ and $S_{\#}\subset\cE_{\#}$, there are
no points $\tilde{x}$ having diagram $(S,S_{\#})$.  Here we determine
exactly which diagrams can arise.

\begin{proposition}\label{prop:triangle-rules}
Let $\Gamma$ be the diagram of a point $\tilde{x}\in\tilde{X}$. Then
the following conditions must hold for $\Gamma$.
\begin{enumerate}{}
\item[(i)] For each of the connected components in
Figure~\ref{fig:fiber} at least one edge must be unmarked (i.e., not bold).
\item[(ii)] Up to symmetry, the set of bold edges in any pair of related
triangles must be one of the $5$ options in Figure \ref{fig:vertrels}.
\begin{figure}[ht]
\begin{center}
\includegraphics[scale = .3]{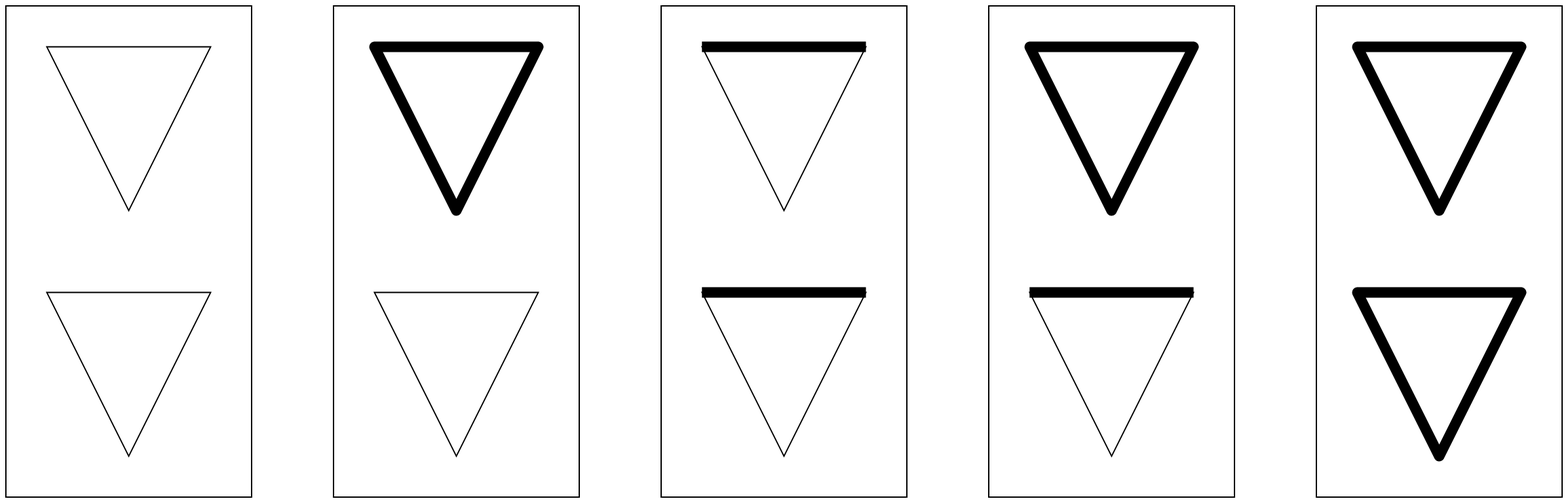}
\end{center}
\caption{\label{fig:vertrels}}
\end{figure}
\end{enumerate}
\end{proposition}

\begin{proof}{}
These are all consequences of Proposition~\ref{prop:local-diagram} and
the linear and quadric equations in Theorem~\ref{thm:local-equations}.
\end{proof}

\begin{definition}\hspace{1in}
\begin{itemize}
\item A diagram $(S,S_{\#})$ with $S\subset\cE$ and
$S_{\#}\subset\cE_{\#}$ is {\em admissible} if it satisfies (i) and
(ii) in the statement of Proposition \ref{prop:triangle-rules}.
\item An admissible diagram $(S,S_{\#})$ is a {\em shifting} diagram
if $S_{\#}=\emptyset$ and the number of vertices in each hypersimplex
$\Delta_k$ is either $1$ or $\binom{4}{k}$ after collapsing all of the
edges in $S$.
\item An admissible diagram $(S,S_{\#})$ is {\em split} if the number
of vertices in each hypersimplex $\Delta_k$ is $\geq 2$ after
collapsing all of the edges in $S$.
\end{itemize}
\end{definition}

The classification of shifting diagrams and split diagrams is a
pleasant combinatorial excercise using
Proposition~\ref{prop:triangle-rules}.  The results are given
Tables~\ref{tab:shifting} and \ref{tab:split}, listed by combinatorial
type.  In the shifting diagrams, we also use the shorthand $\bullet$
(respectively $\circ $) to indicate that an affine hypersimplex
component is collapsed (resp., not collapsed).  Note that
one diagram, $X_{\emptyset}$, is both shifting and split; this is the
only diagram having this property.

\begin{table}[ht]\caption{\label{tab:shifting} The shifting diagrams.}
\begin{tabular}{||c|c|c|c|c||c||}\hline
Type &  $S$ & $S_{\#}$ & $\cod$ & \# \\
\hline\hline
$X_{\emptyset}\quad (\circ \circ \circ) $ &\includegraphics[scale = .3]{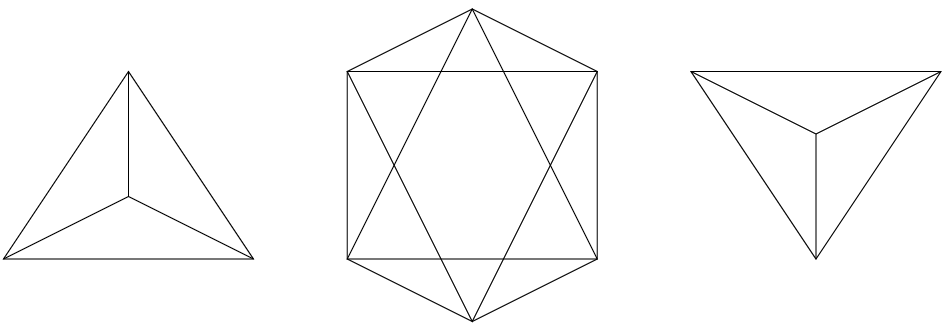} &
\includegraphics[scale = .27]{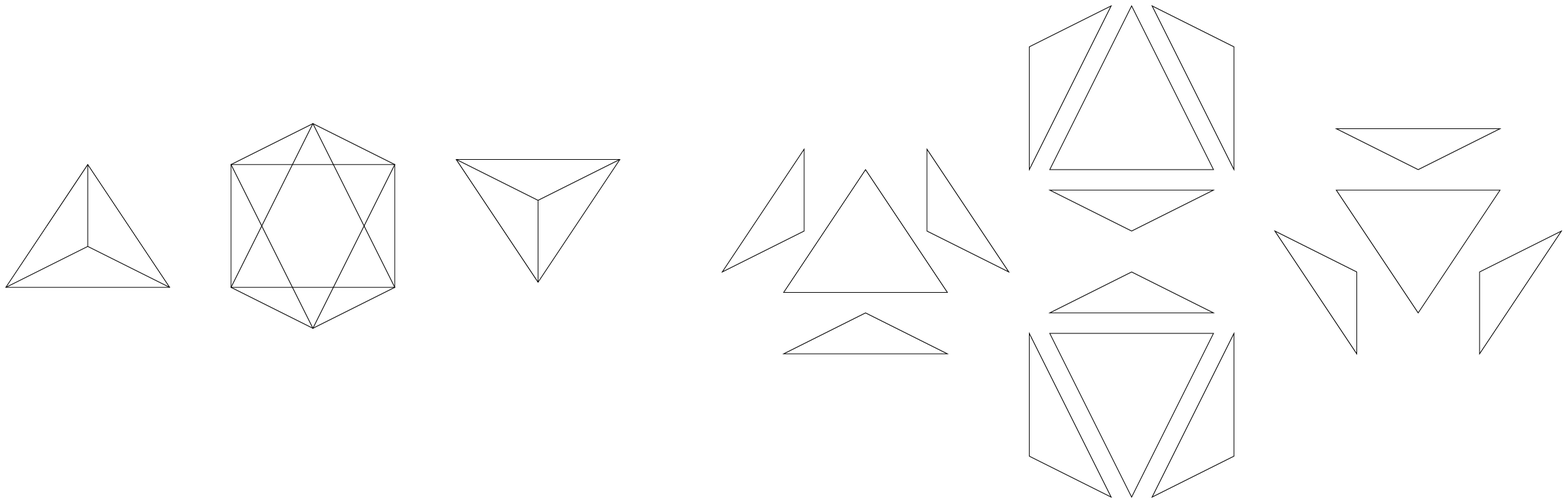} & $0$ &$1$ \\
\hline\hline
$A\quad (\bullet\circ \circ )$ & \includegraphics[scale = .3]{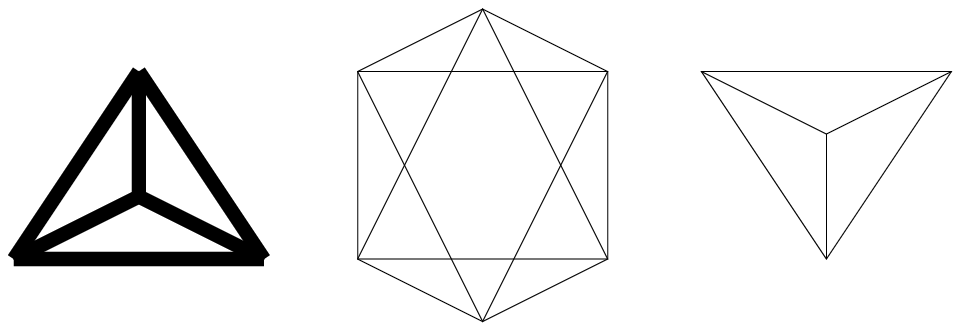} &
\includegraphics[scale = .27]{fiber-gen.eps}
& $1$ & $1$\\
\hline
$B\quad (\circ \bullet \circ )$ & \includegraphics[scale = .3]{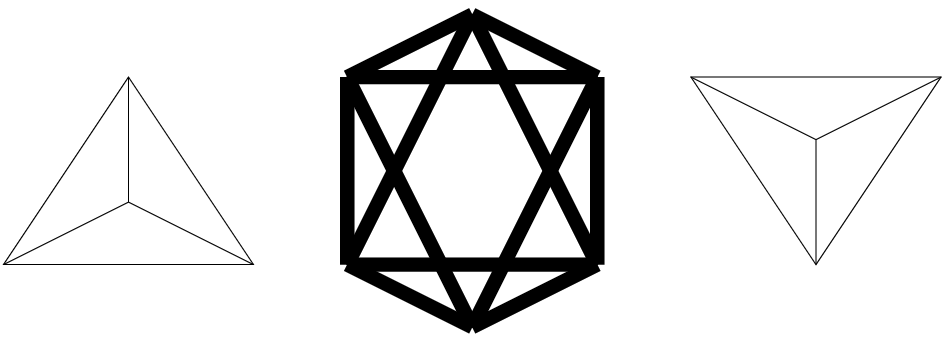} &
\includegraphics[scale = .27]{fiber-gen.eps}
& $1$ & $1$\\
\hline
$A^{*}\quad (\circ \circ \bullet)$ & \includegraphics[scale = .3]{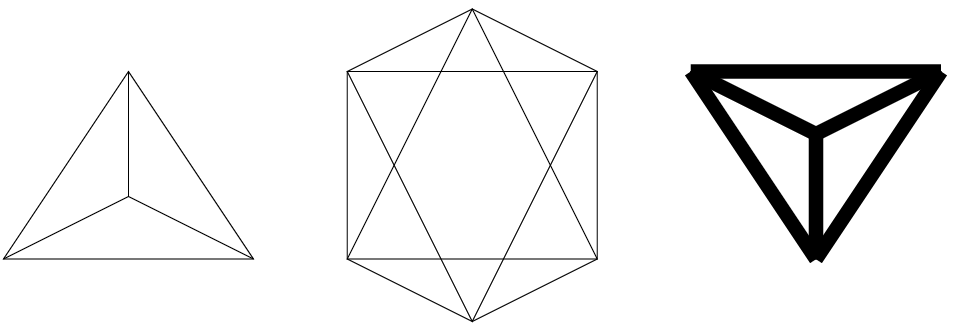} &
\includegraphics[scale = .27]{fiber-gen.eps}
& $1$ & $1$\\
\hline\hline
$AB\quad (\bullet\bullet\circ)$ & \includegraphics[scale = .3]{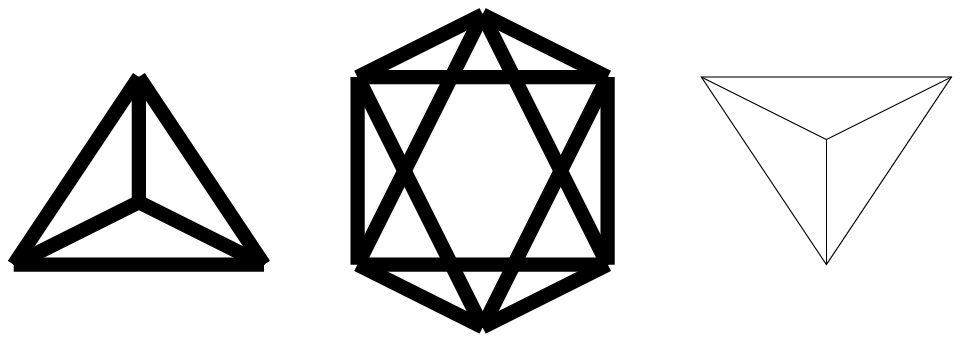} &
\includegraphics[scale = .27]{fiber-gen.eps}
& $2$ & $1$\\
\hline
$AA^{*}\quad (\bullet\circ \bullet)$ & \includegraphics[scale = .3]{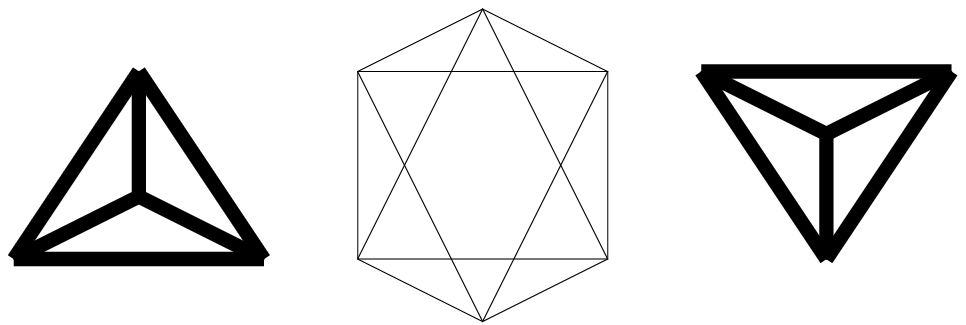} &
\includegraphics[scale = .27]{fiber-gen.eps}
& $2$ & $1$\\
\hline
$BA^{*}\quad (\circ \bullet\bullet)$ & \includegraphics[scale = .3]{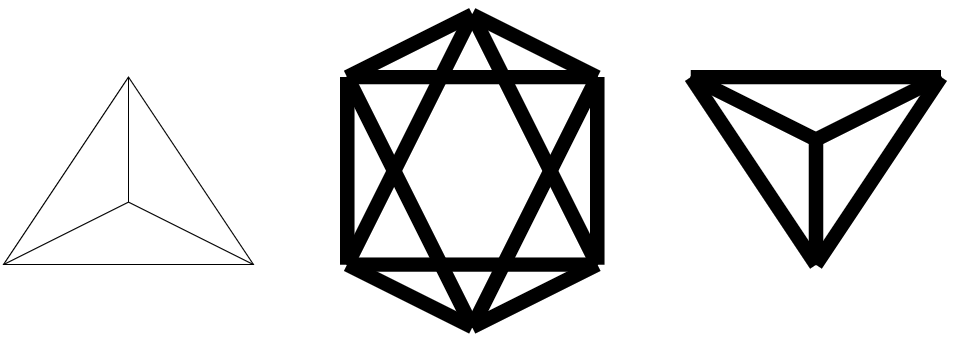} &
\includegraphics[scale = .27]{fiber-gen.eps}
& $2$ & $1$\\
\hline\hline
$ABA^{*}\quad(\bullet\bullet\bullet)$ & \includegraphics[scale = .3]{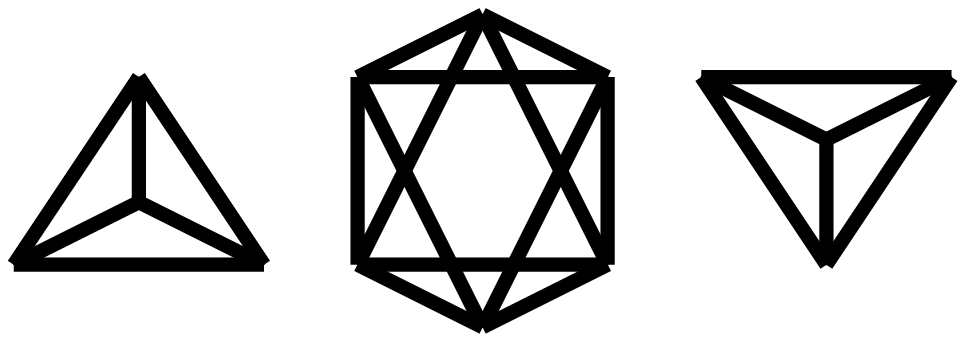} &
\includegraphics[scale = .27]{fiber-gen.eps}
& $3$ &
$1$\\
\hline\hline
\end{tabular}
\end{table}

\begin{table}[t]\caption{\label{tab:split} The split diagrams.}
\begin{tabular}{||c|c|c|c|c||c||} \hline
Type &  $S$ &  $S_{\#}$ &
$\cod$ & \# \\
\hline\hline
$X_{\varnothing}$ &\includegraphics[scale = .3]{diagram-gen.eps} &
 \includegraphics[scale = .27]{fiber-gen.eps} &
$0$ & $1$\\
\hline\hline
$C$ & \includegraphics[scale = .3]{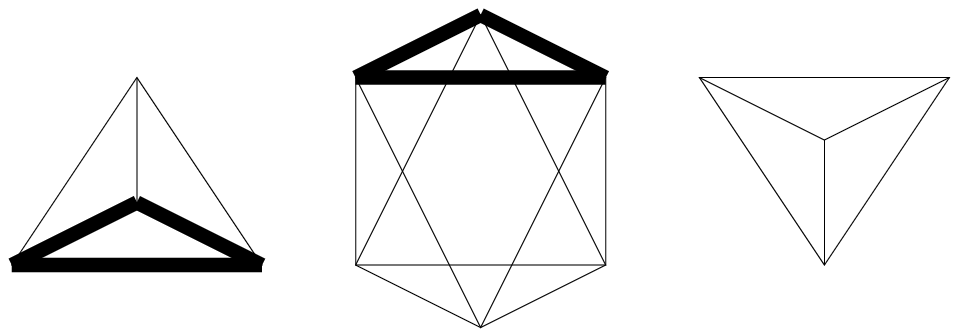} &
\includegraphics[scale = .27]{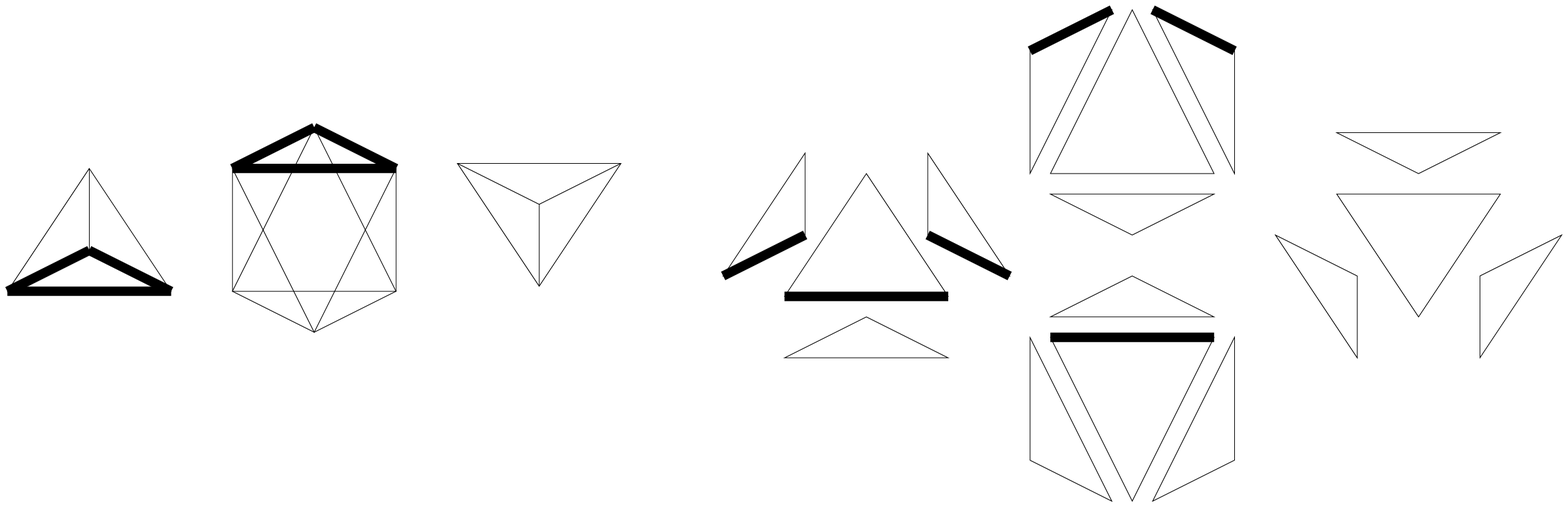} & $1$ &
$4$\\
\hline
$C^*$ & \includegraphics[scale = .3]{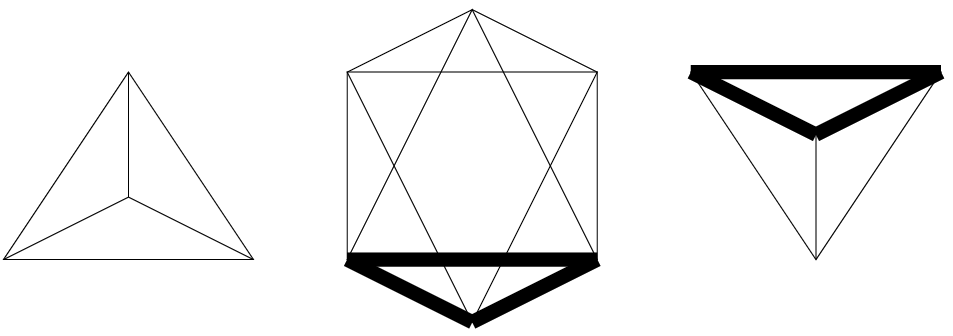} & \includegraphics[scale = .27]{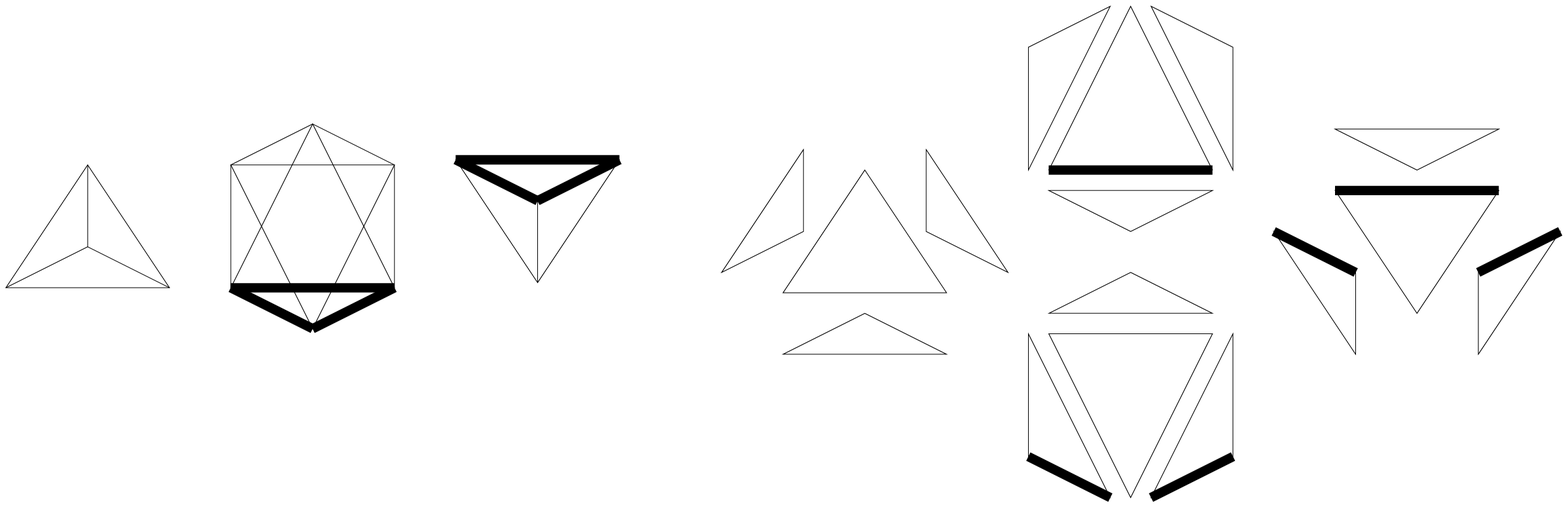} &  $1$ & $4$\\
\hline
$D$ & \includegraphics[scale = .3]{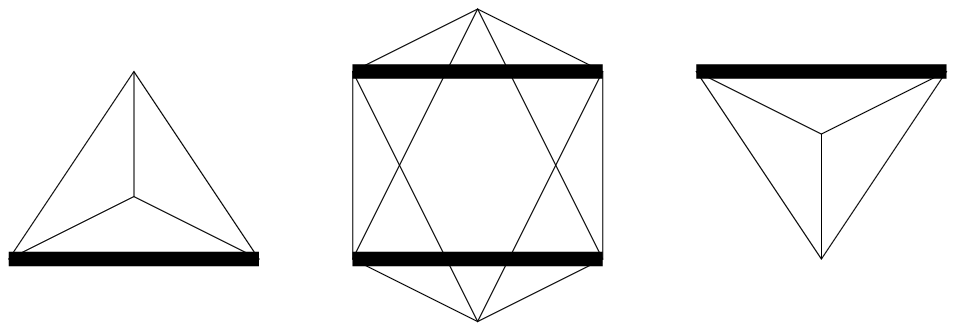} &  \includegraphics[scale = .27]{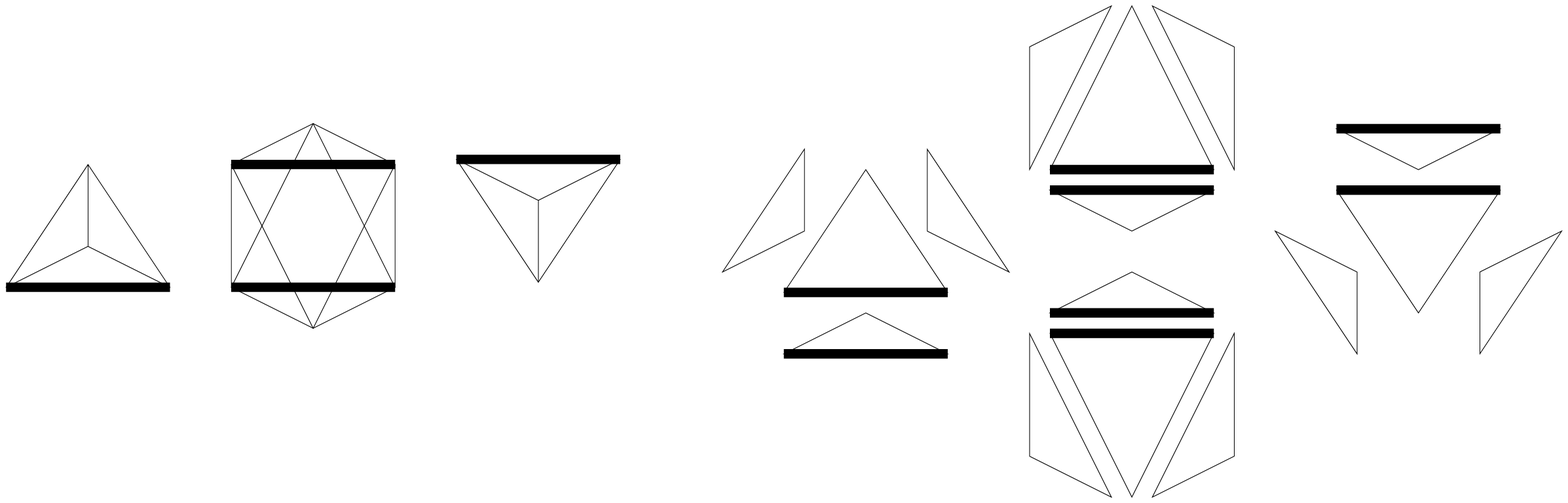} &$1$& $6$\\
\hline
$E$ & \includegraphics[scale = .3]{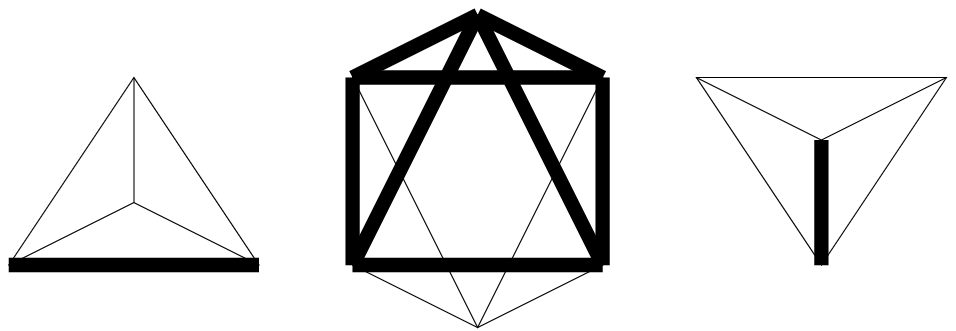} &
\includegraphics[scale = .27]{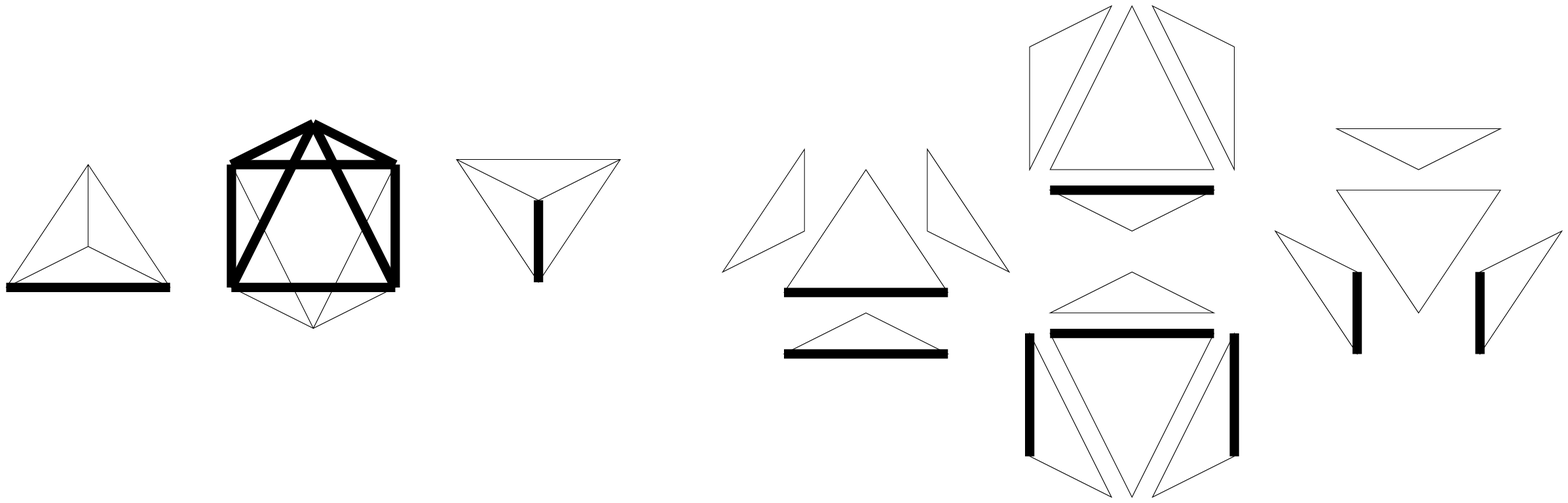} & $1$ & $6$\\
\hline\hline
$CC_{nop}^*$ & \includegraphics[scale =
.3]{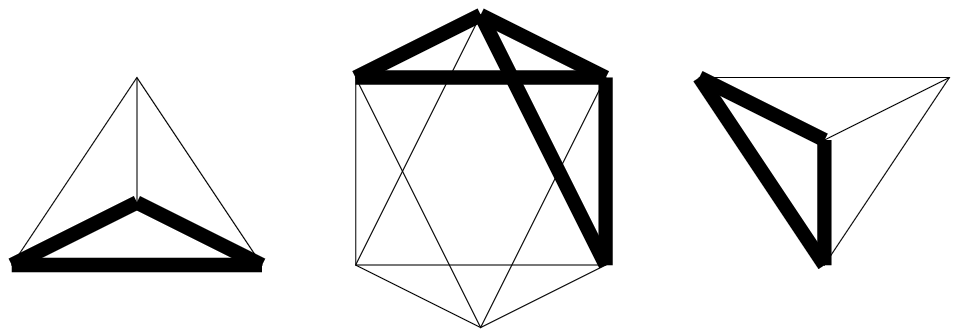}
& \includegraphics[scale = .27]{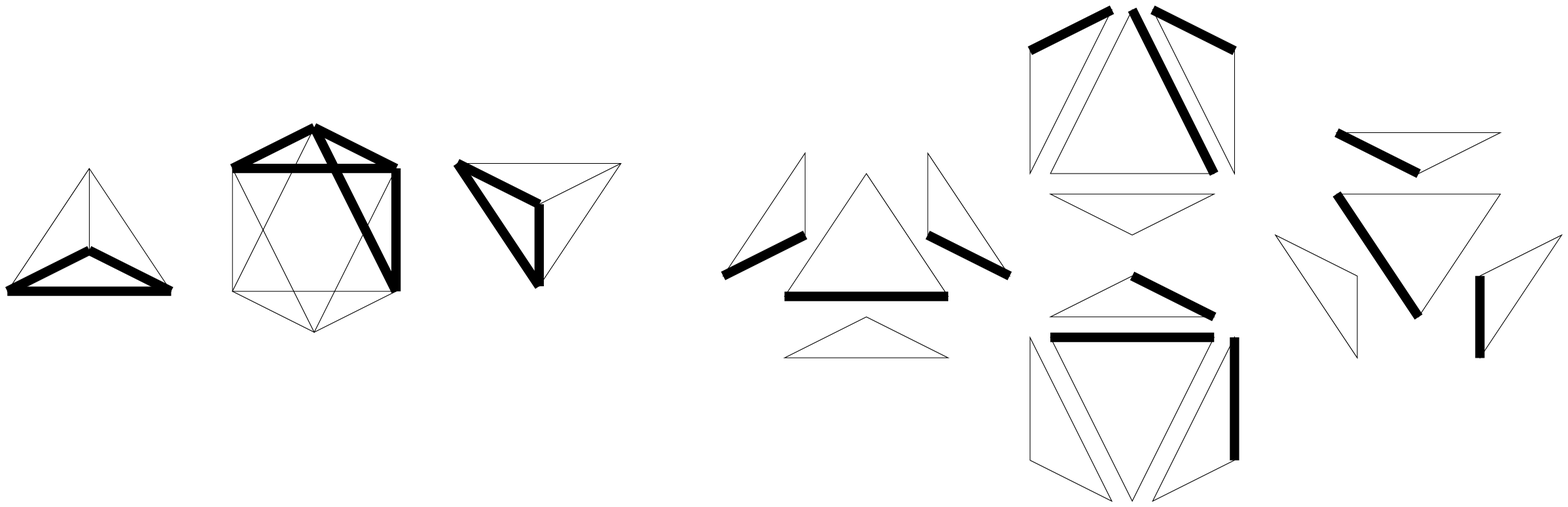}& $2$ &$12$\\
\hline
$CC_{op}^*$ & \includegraphics[scale =
.3]{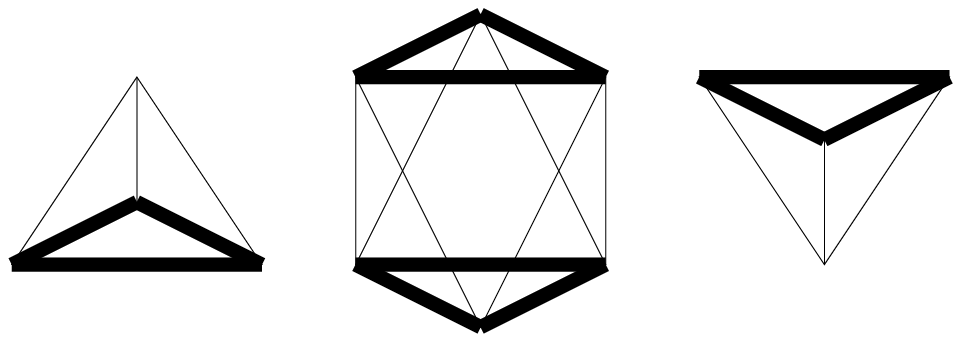} &
\includegraphics[scale = .27]{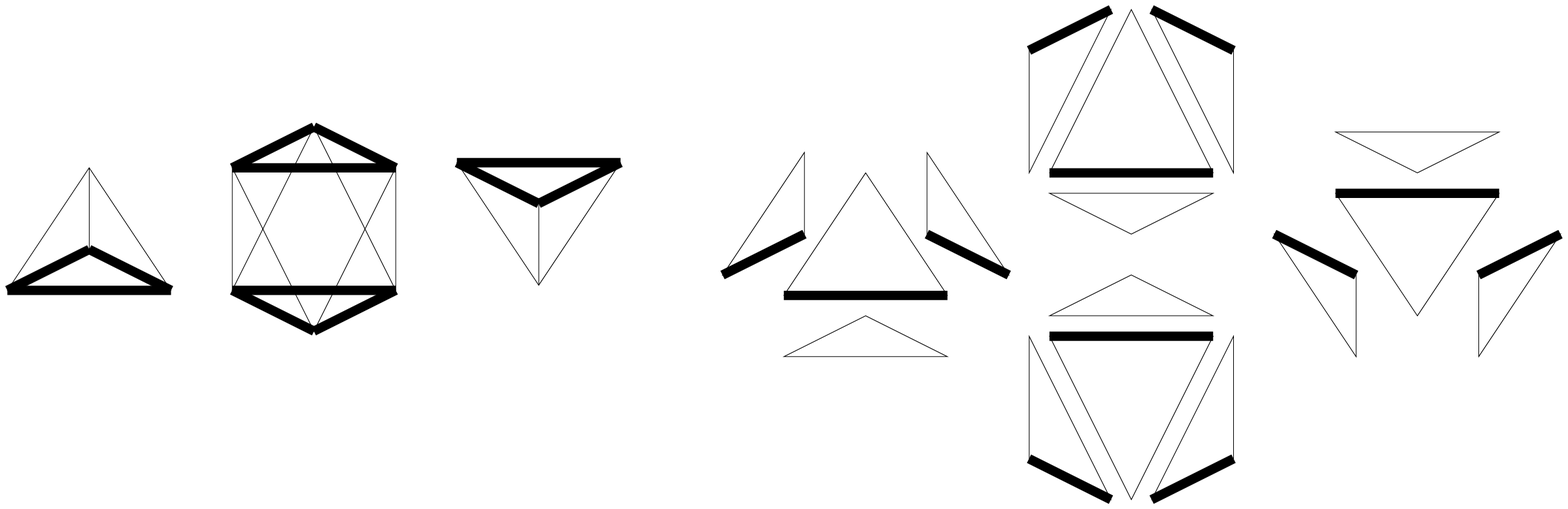} & $2$ &$4$\\
\hline
$CD$ & \includegraphics[scale = .3]{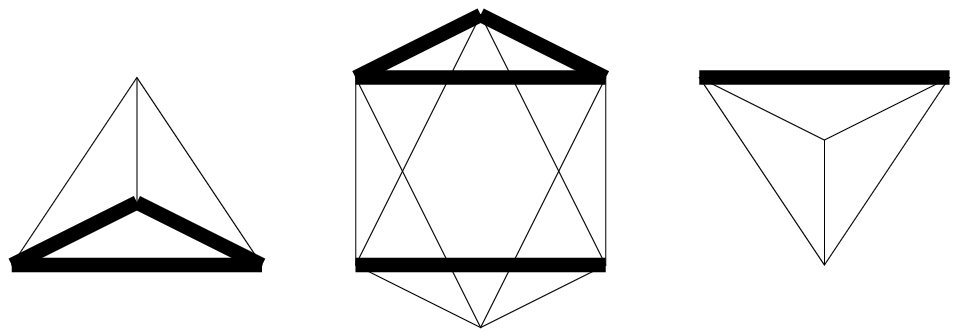} &
 \includegraphics[scale = .27]{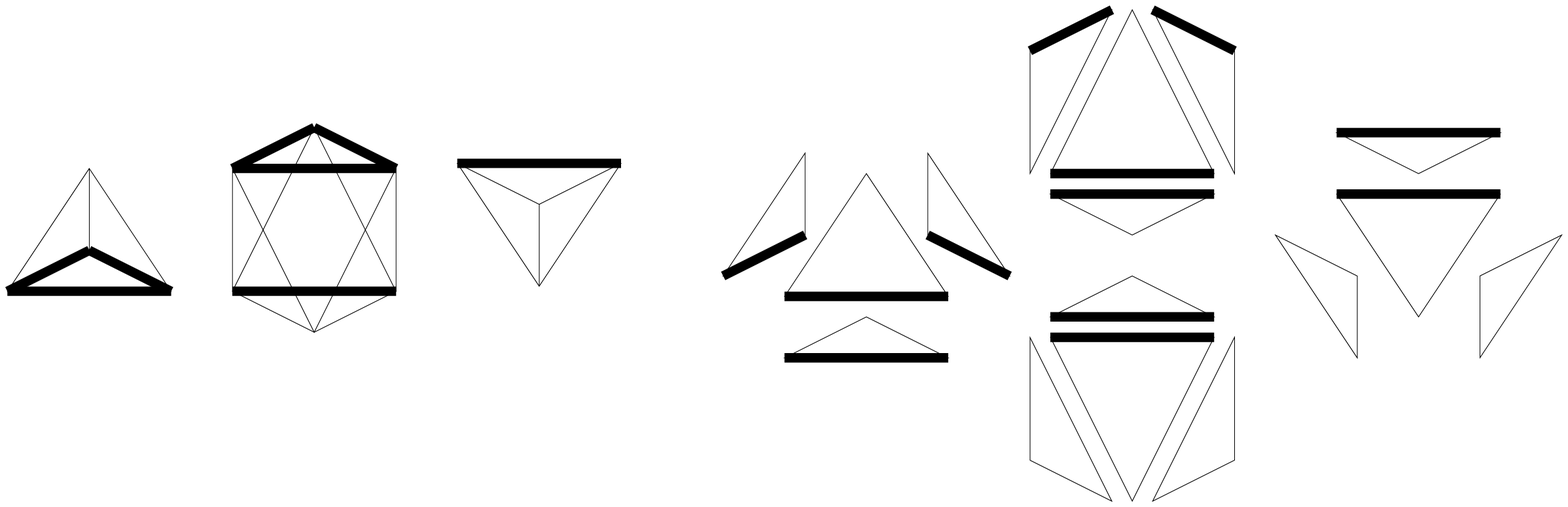} & $2$  &
$12$\\
\hline
$C^*D$ & \includegraphics[scale = .3]{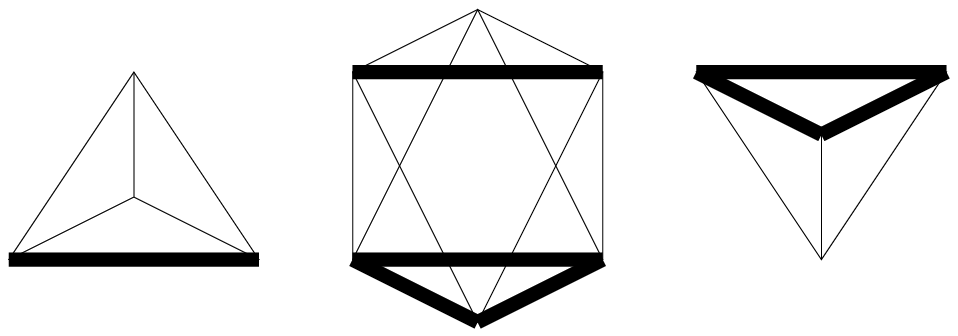} &
 \includegraphics[scale = .27]{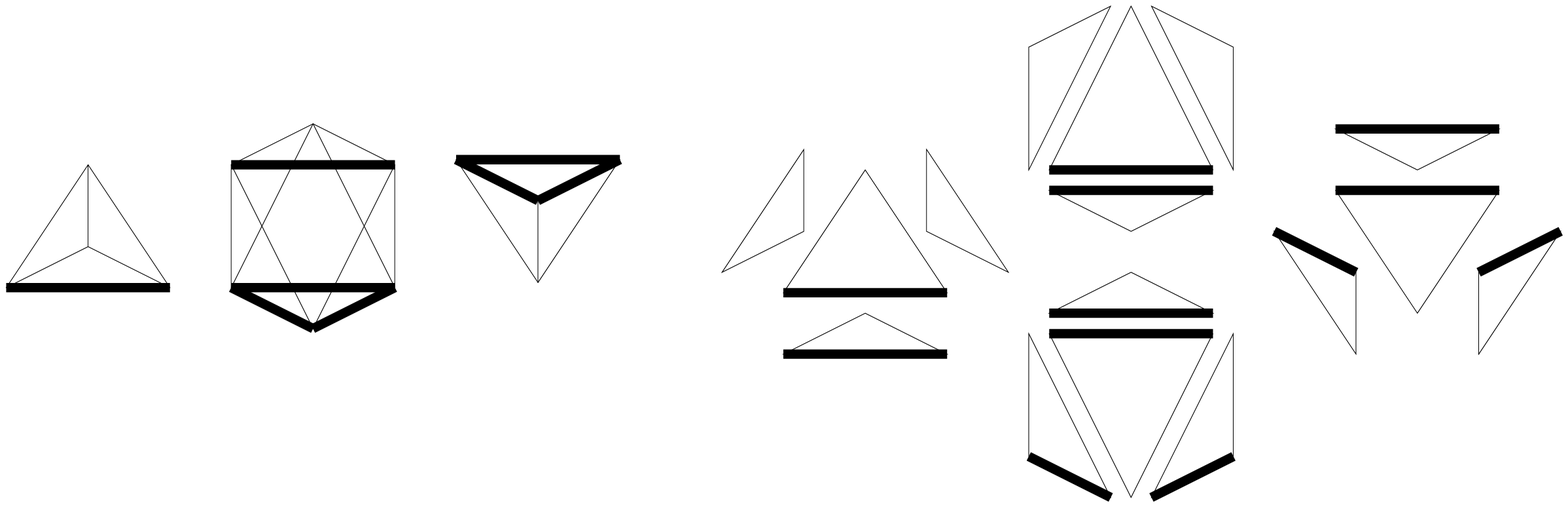} & $2$ &
$12$\\
\hline
\end{tabular}
\end{table}

\begin{table}[ht]
\begin{tabular}{||c|c|c|c|c||c||} \hline
Type &  $S$ &  $S_{\#}$ &
$\cod$ & \# \\
\hline\hline
$CE$ & \includegraphics[scale = .3]{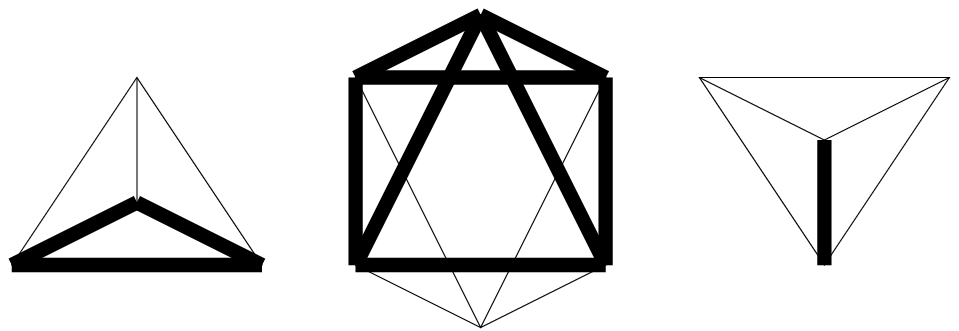} &
 \includegraphics[scale = .27]{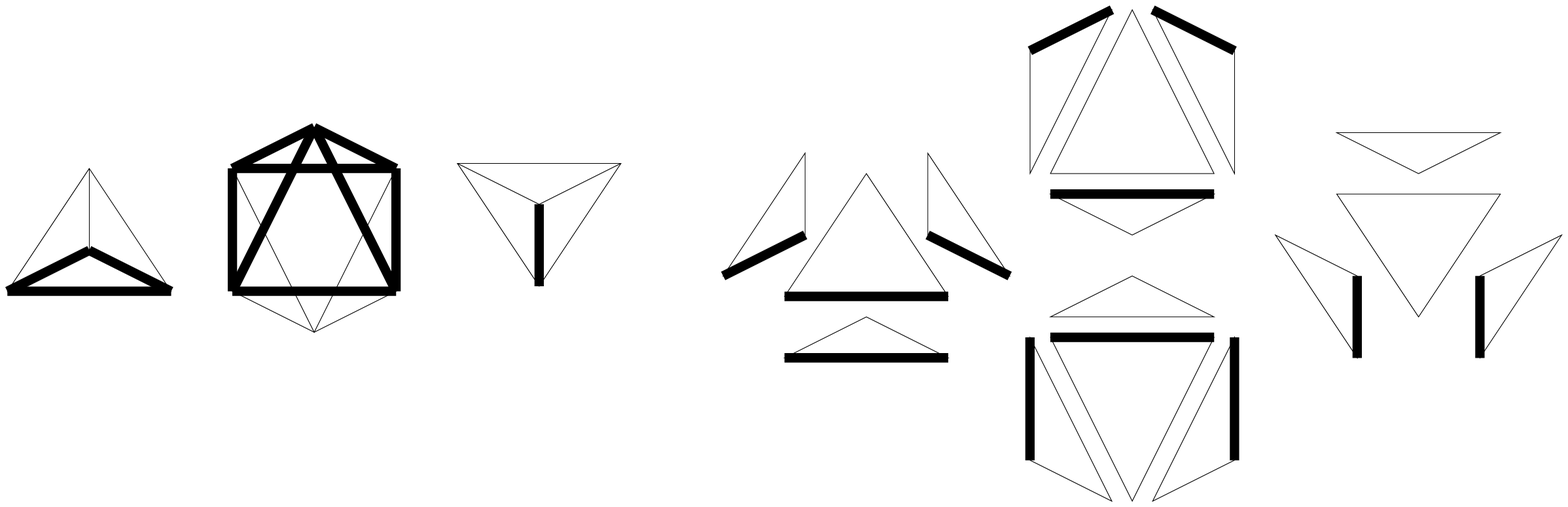} & $2$ &
$12$\\
\hline
$C^*E$ & \includegraphics[scale = .3]{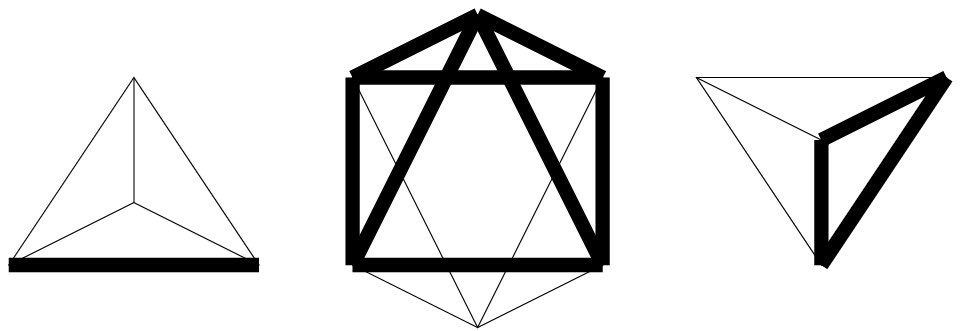} &
\includegraphics[scale = .27]{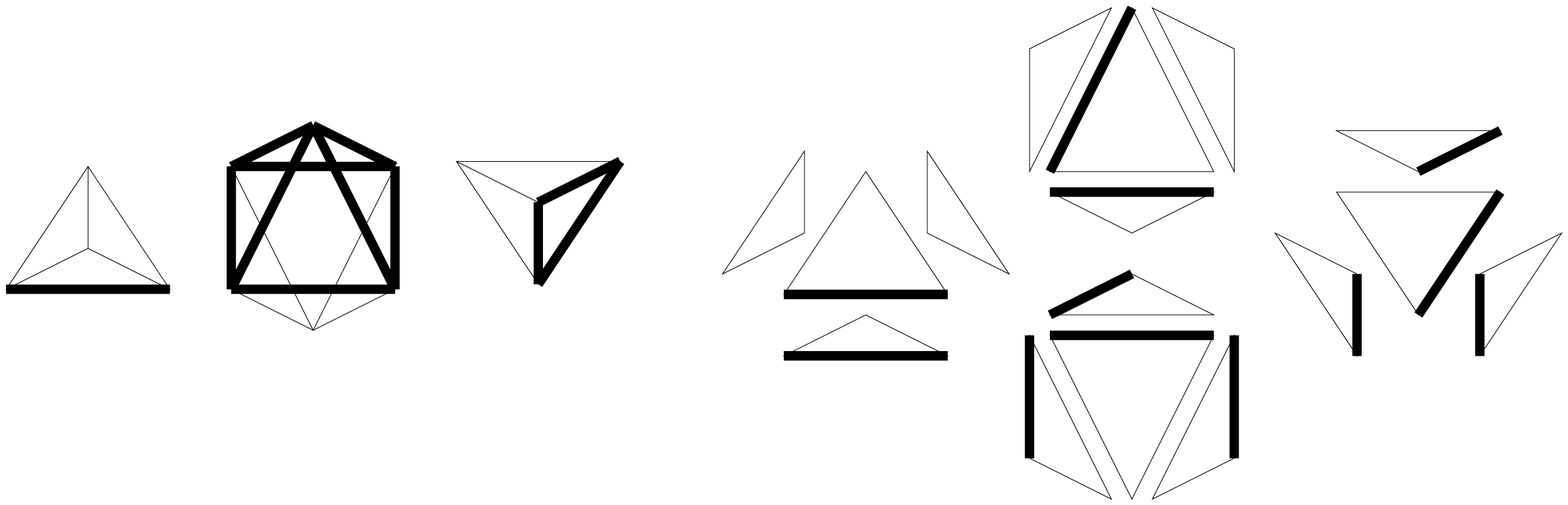} & $2$ &
$12$\\
\hline
$DD_{op}$ & \includegraphics[scale = .3]{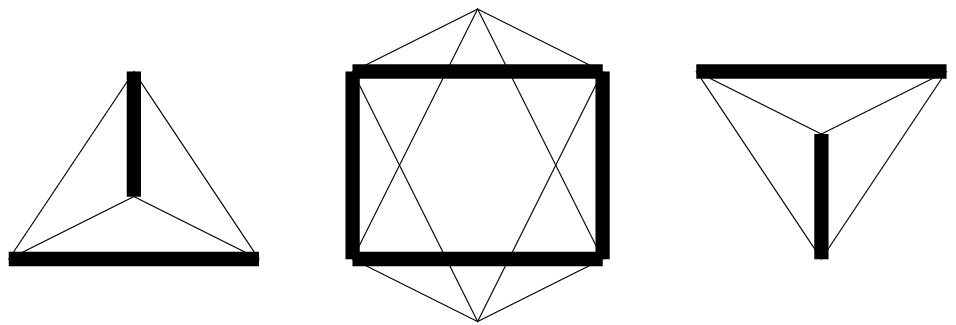} &
\includegraphics[scale = .27]{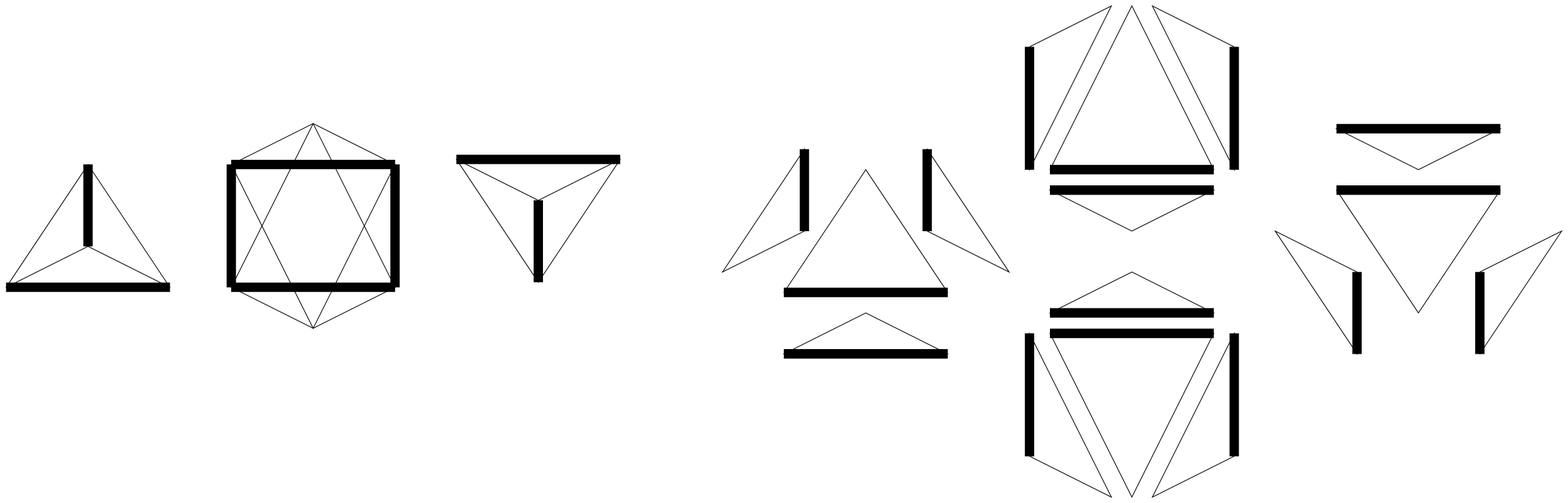} & $2$ &
$3$\\
\hline
$DE$ & \includegraphics[scale = .3]{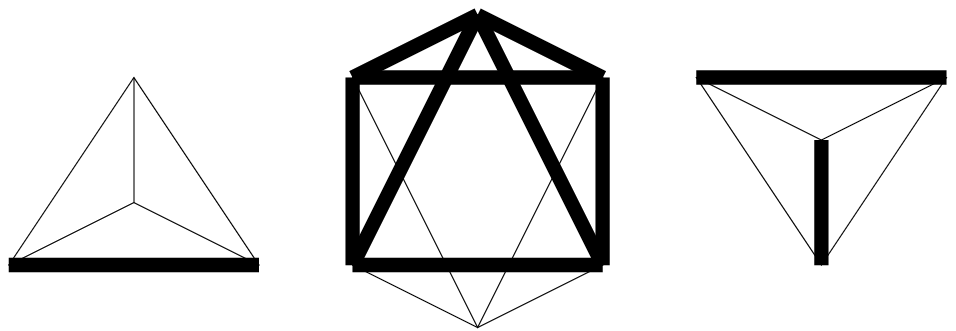} &
\includegraphics[scale = .27]{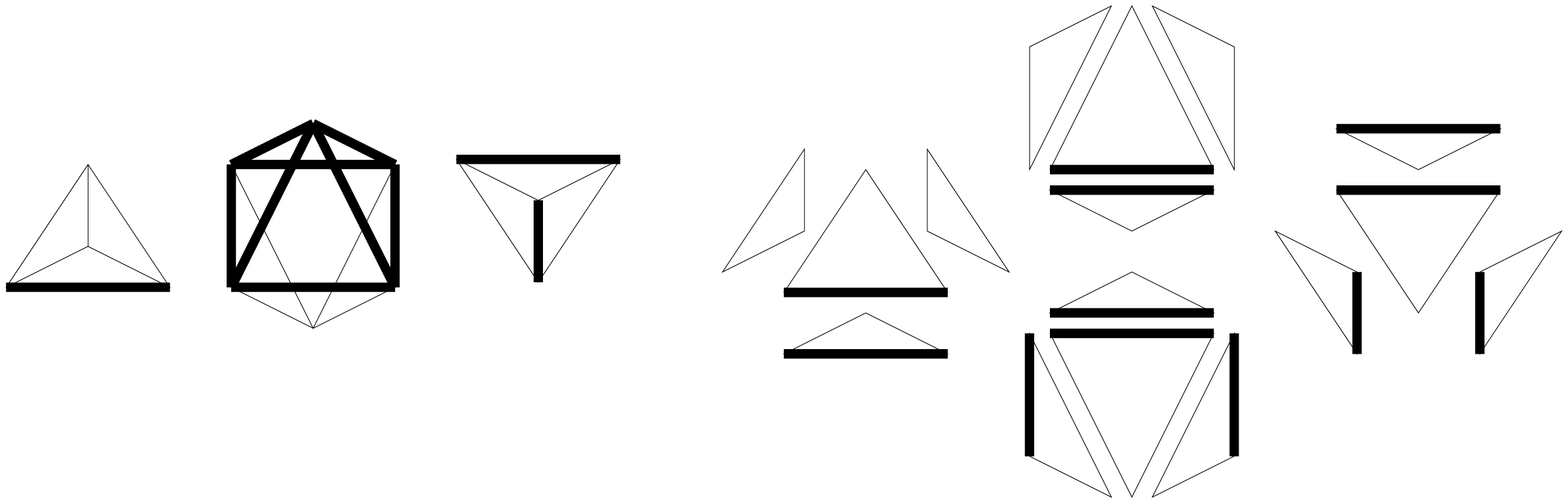} & $2$ &
$6$\\
\hline
$DE^*$ & \includegraphics[scale = .3]{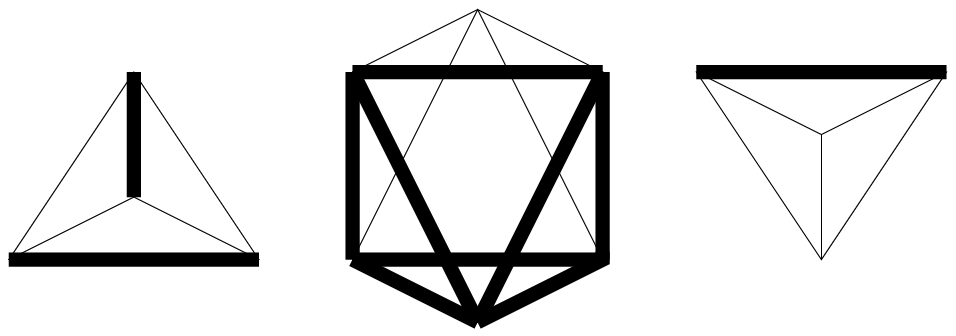} &
\includegraphics[scale = .27]{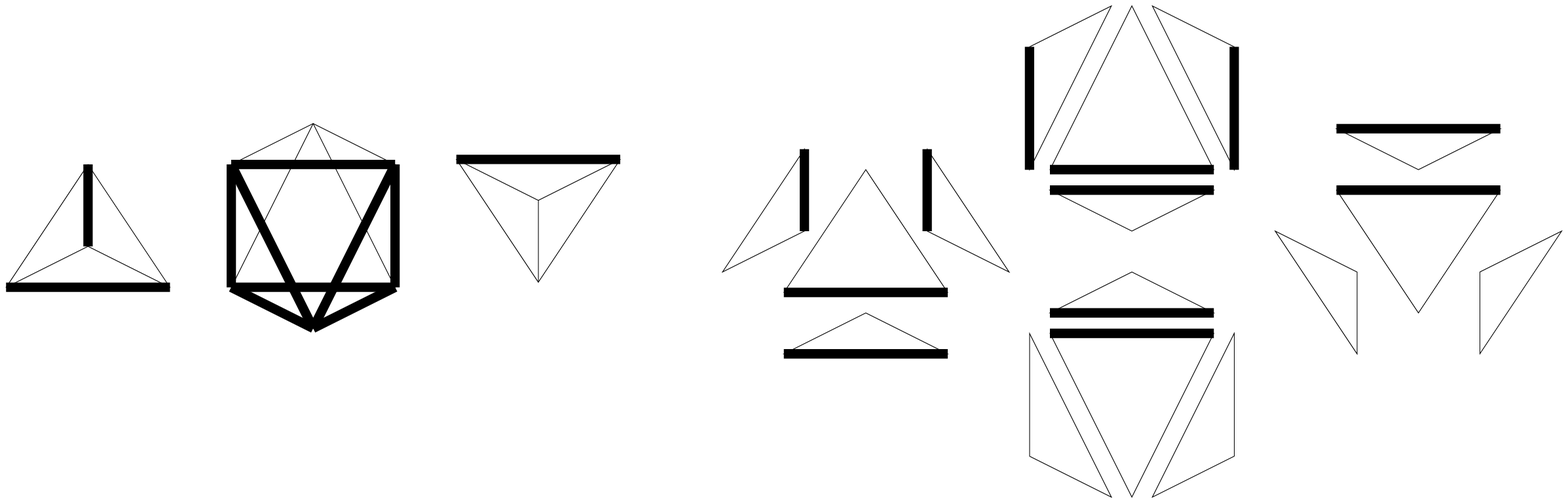} & $2$ &
$6$\\
\hline\hline
$CC_{op}^*D$ & \includegraphics[scale =
.3]{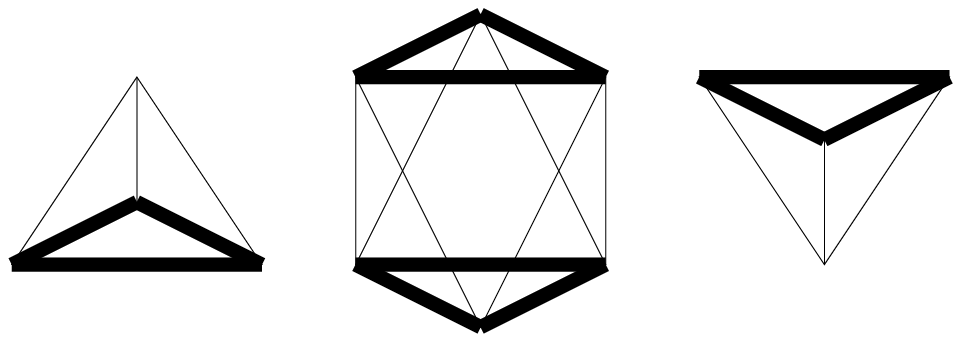} & \includegraphics[scale =
.27]{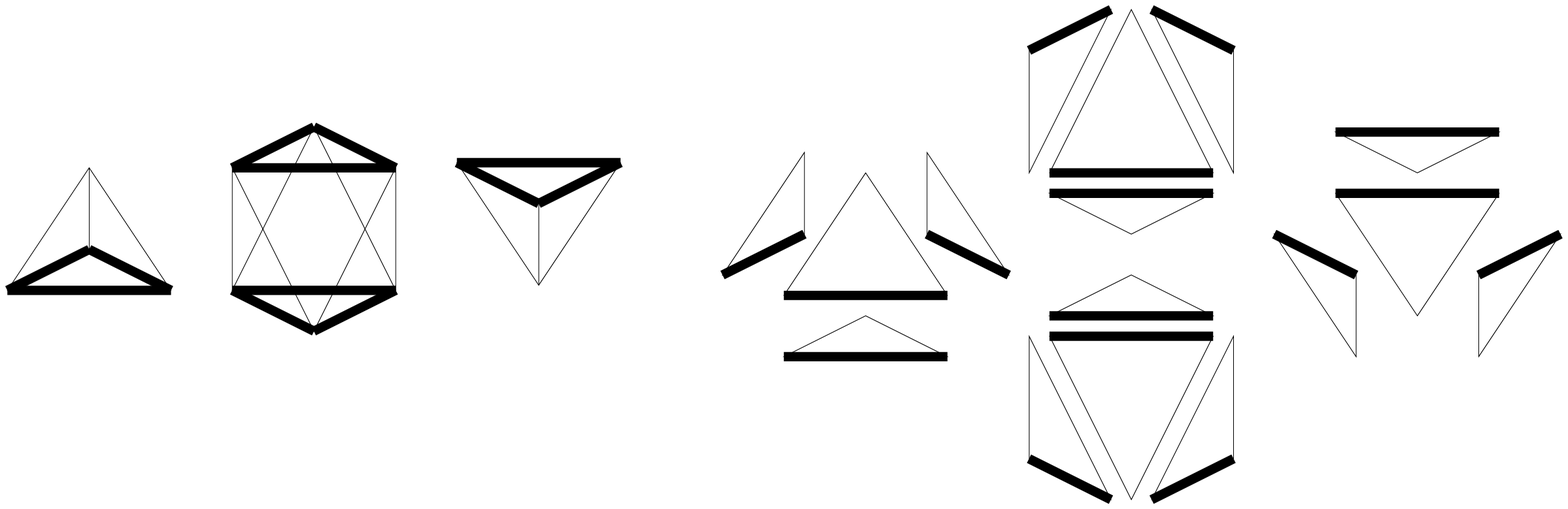} & $3$ & $12$\\
\hline
$CC_{nop}^*D$ & \includegraphics[scale = .3]{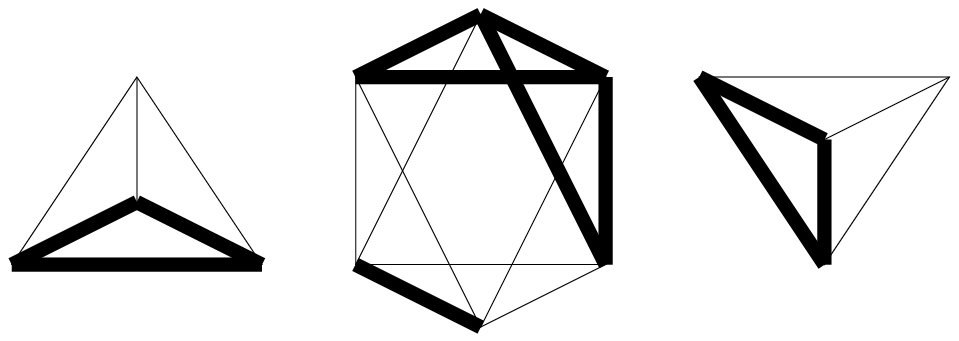}
& \includegraphics[scale = .27]{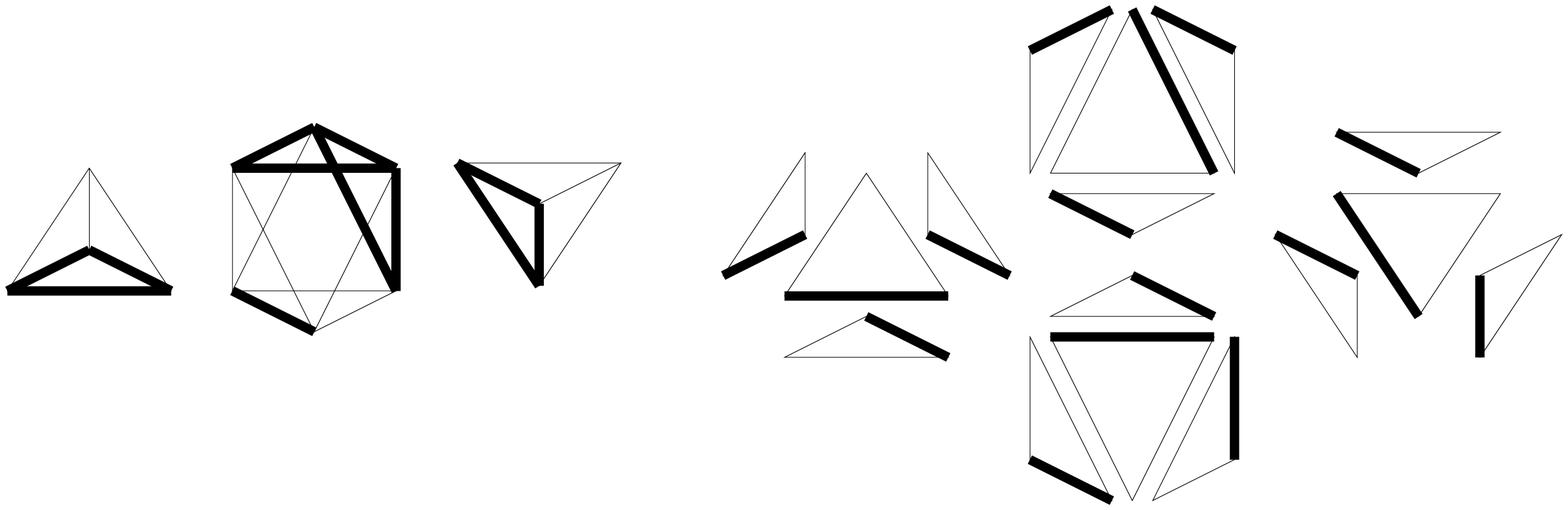} & $3$ & $12$\\
\hline
$CC^*E$ & \includegraphics[scale = .3]{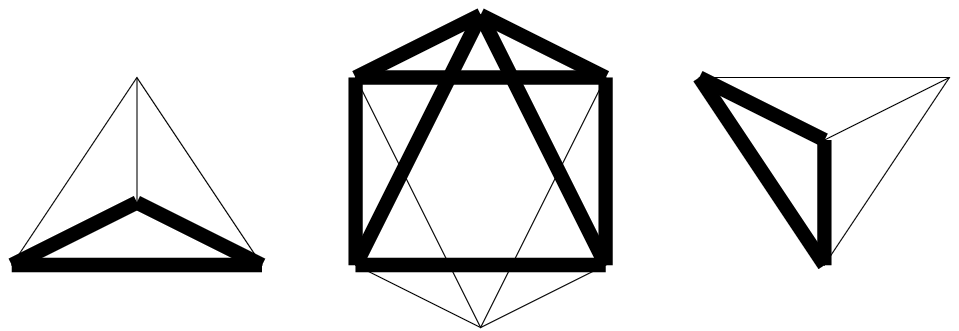} &
\includegraphics[scale = .27]{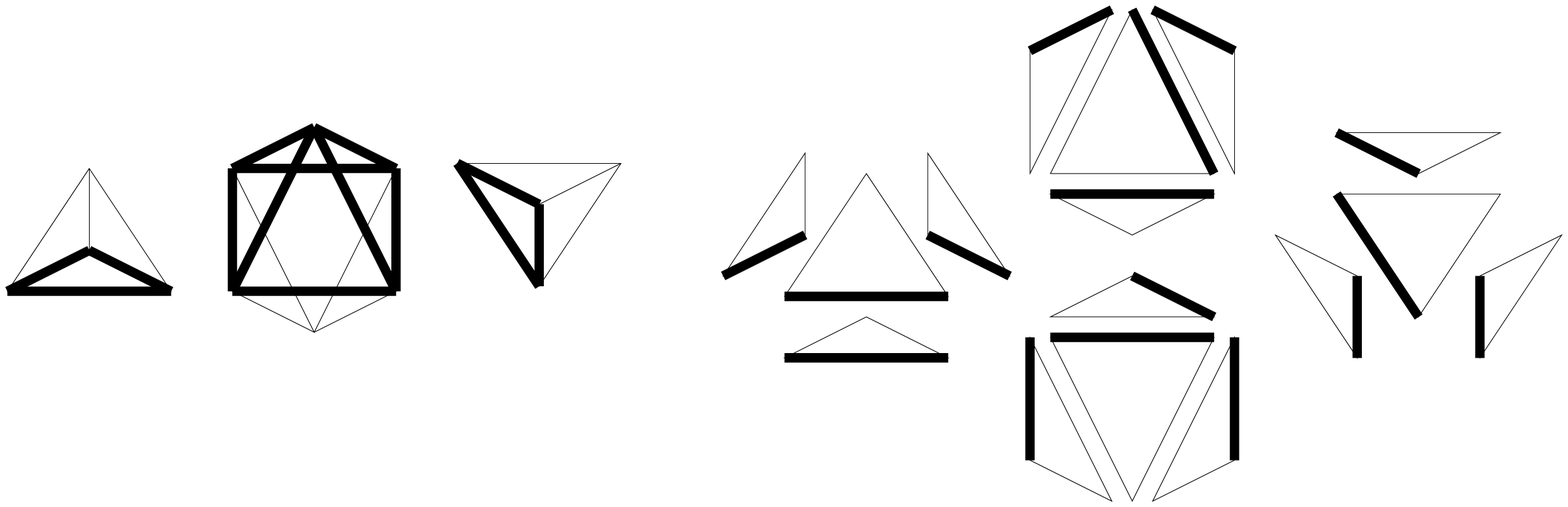} & $3$ &
$24$\\
\hline
$CDE$ & \includegraphics[scale = .3]{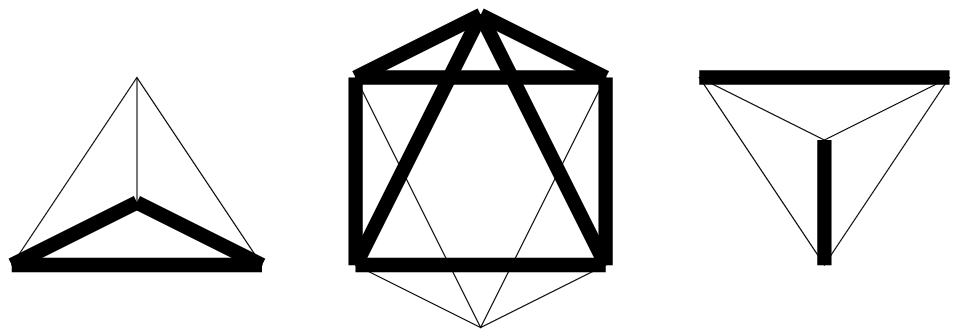} &
\includegraphics[scale = .27]{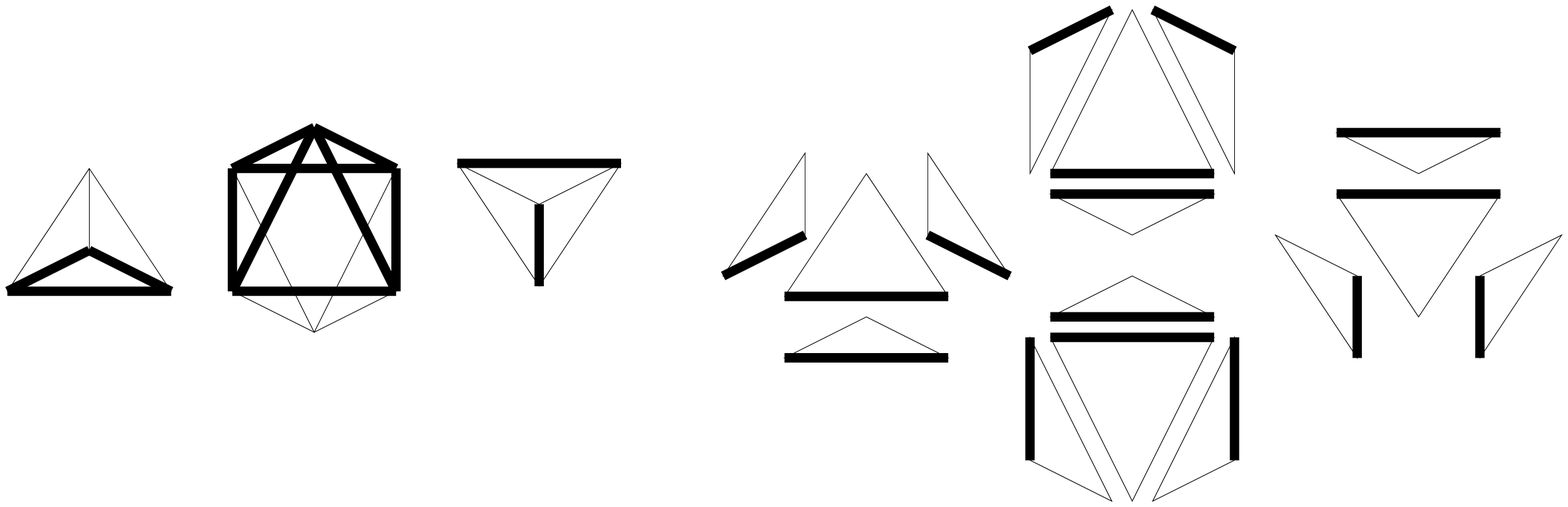} & $3$ &
$12$\\
\hline
\end{tabular}
\end{table}

\begin{table}[ht]
\begin{tabular}{||c|c|c|c|c||c||} \hline
Type &  $S$ &  $S_{\#}$ &
$\cod$ & \# \\
\hline\hline
$C^*DE$ & \includegraphics[scale = .3]{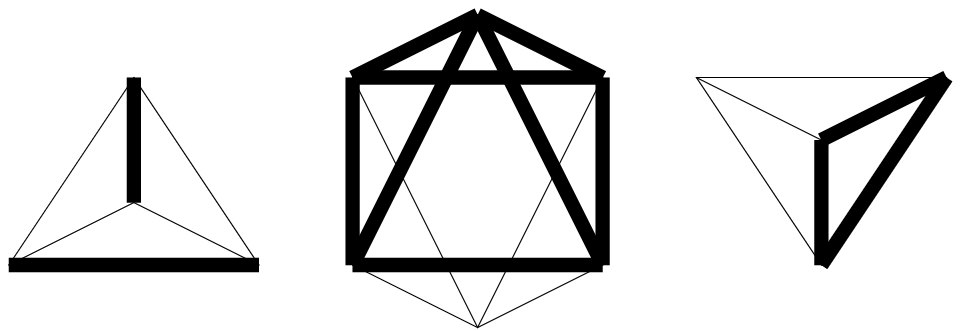} &
\includegraphics[scale = .27]{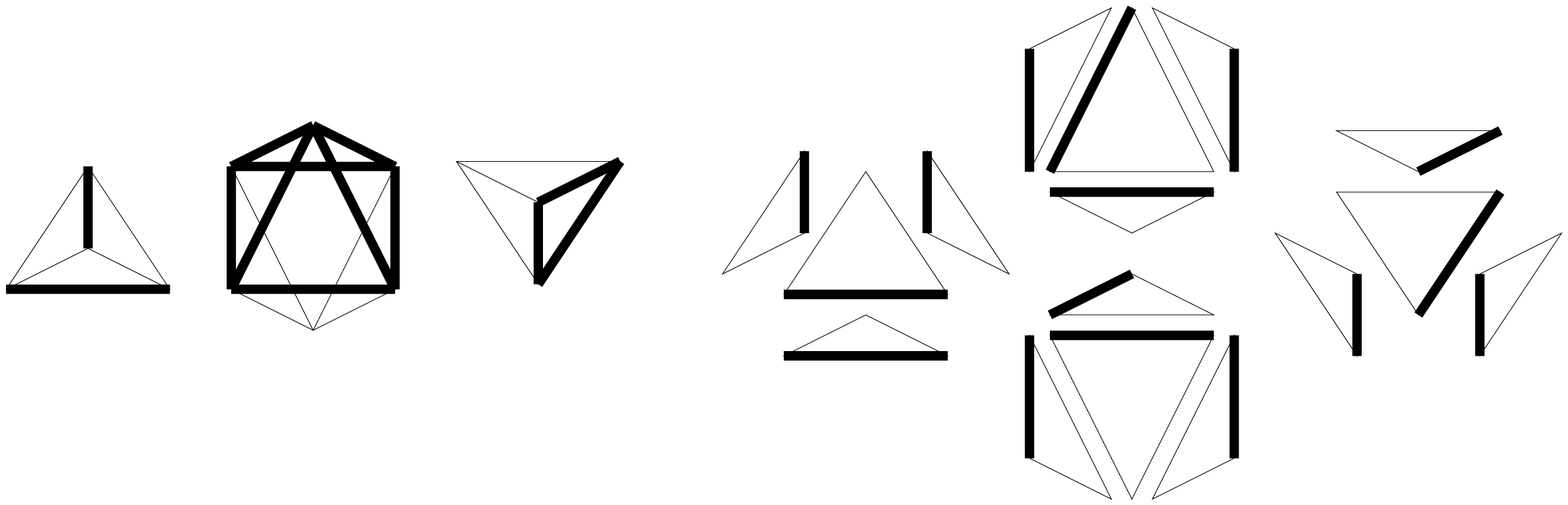} & $3$ &
$12$\\
\hline
$DD_{op}E$ & \includegraphics[scale = .3]{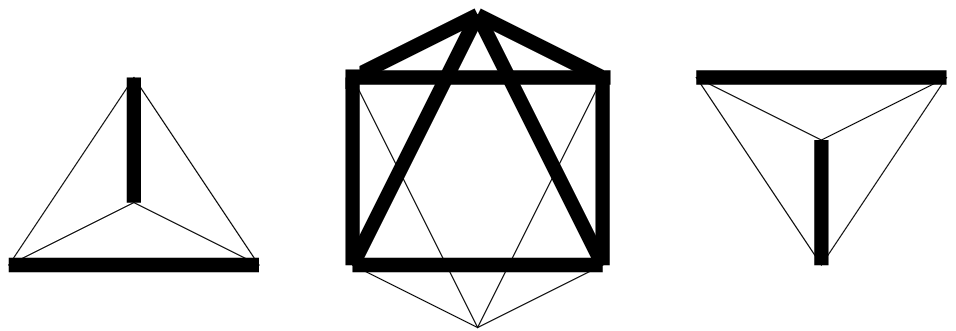} &
\includegraphics[scale = .27]{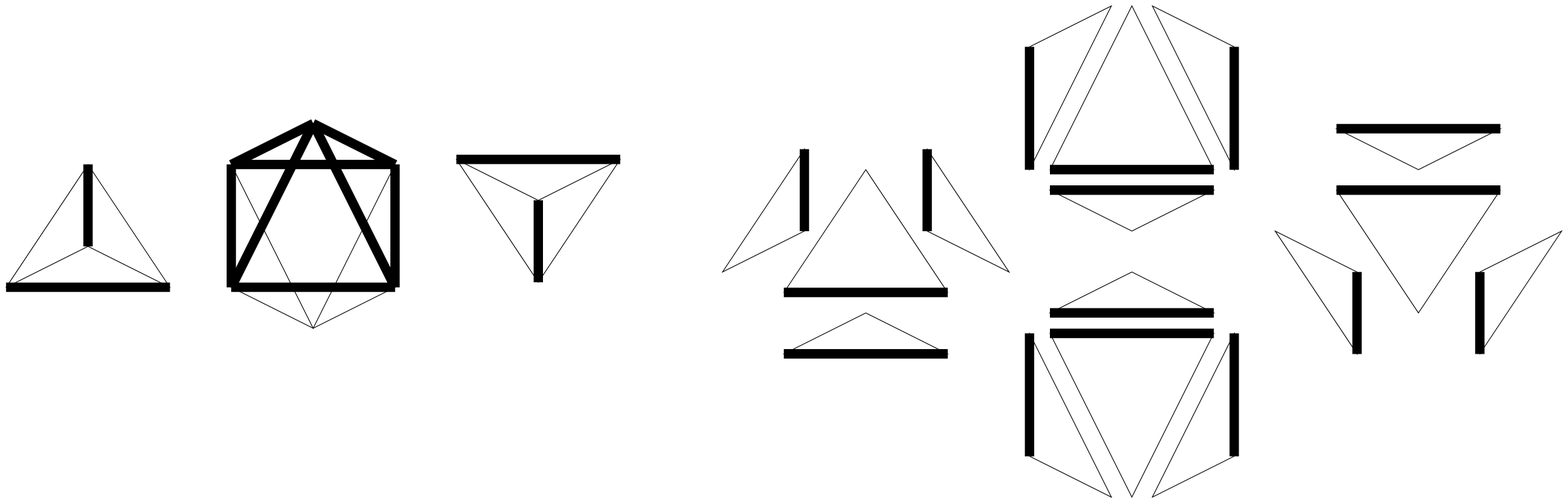} & $3$ &
$6$\\
\hline\hline
\end{tabular}
\end{table}
\clearpage

The classification of admissible diagrams can be formulated in
terms of shifting and split diagram as follows:

\begin{proposition}\label{prop:diagram-classification}
A diagram $\Gamma$ is admissible if and only if it can be written
in the form $\Gamma=(S_{\spl}\cup S_{\shi},S_{\#})$ where
$(S_{\shi},\emptyset)$ is a shifting diagram (called the
``shifting part'' of $\Gamma$) and $(S_{\spl},S_{\#})$ is a split
diagram (called the ``split part'' of $\Gamma$).  Moreover, this
representation is unique.
\end{proposition}

\begin{proof}{}
The diagram $(S_{\spl}\cup S_{\shi},S_{\#})$ is obtained from the
split diagram $(S_{\spl},S_{\#})$ by placing the fully marked
hypersimplices in $S_{\shi}$ into the corresponding location in the diagram
for $S_{\spl}$, and leaving the remaining components of $S_{\spl}$
and $S_{\#}$ unchanged.  The resulting marked diagram satisfies (i)
and (ii) of Proposition~\ref{prop:triangle-rules}, and hence is
admissible.

Conversely, suppose $\Gamma=(S,S_{\#})$ is admissible.  The rules (i)
and (ii) imply that the only possibilities for $S_{\#}$ are those
appearing in Table~\ref{tab:split}, and each such $S_{\#}$ appears
exactly once.  Letting $(S_{\spl},S_{\#})$ be the unique
corresponding split diagram, we see that the same rules imply that $S$
must contain all of the edges in $S_{\spl}$, and that adding any
other edge to $S_{\spl}$ from the hypersimplex $\Delta_i$ forces us
to add all of the edges in $\Delta_i$.  Thus, $S=S_{\spl}\cup
S_{\shi}$ for exactly one of the $8$ shifting diagrams
$(S_{\shi},\emptyset)$.

The uniqueness of the decomposition into split and shifting parts is
clear, and this completes the proof.
\end{proof}

\subsection{The stratification}

We now want to use admissible diagrams to construct our
stratification.  First we put a partial order on the set of all
diagrams, by putting $(S,S_{\#})<(S',S_{\#}')$ if and only if
$S\supset S'$ and $S_{\#}\supset S_{\#}'$.  Next for any diagram
$\Gamma$, let $\tilde{X}_{\Gamma}$ be the subvariety
\[
\tilde{X}_{\Gamma}=\{\tilde{x}\in \tilde{X}\mid\Gamma(\tilde{x})=\Gamma\}.
\]
We will also be interested in the image $X_{\Gamma}$ of this
subvariety under the projection $\tilde{X}\rightarrow X$.  It follows
immediately from the definition of $\Gamma(\tilde{x})$ that
$X_{\Gamma}$ is the locus
\[X_{\Gamma}=\{x\in X\mid S(x)=S\}\]
where $\Gamma=(S,S_{\#})$.
The set of $\tilde{X}_{\Gamma}$, as $\Gamma$ ranges over the poset of
admissible diagrams $\sS$ equipped with the above partial order, will
be our stratification.  We will abuse language and call the
$\tilde{X}_{\Gamma}$ ``strata,'' even though we have not yet verified
that they form the strata of a stratification.

To understand the meaning of the construction of $\tilde{X}_{\Gamma}$,
recall that in \ref{ss:local-equations} we defined an open cover
$\{\tilde{U} (V) \}_{V\in \Flag}$ of $\tilde{X}$.  Moreover, any
$\tilde{U} (V)$ is embedded in a product of $\Flag$ with certain
affine and projective spaces, and the coordinates of the latter spaces
are indexed by elements of $\cE$.  By Proposition~\ref{prop:local-diagram},
$\tilde{X}_{\Gamma}$ is the locus cut out by the vanishing of the coordinates indexed by $\Gamma$
in \emph{any} of the $\tilde{U} (V)$.  It also follows from
Theorem~\ref{thm:local-equations} and
Proposition~\ref{prop:local-diagram}, that the locally trivial
fibrations $\tilde{X}\rightarrow\Flag$ and $X\rightarrow\Flag$
restrict to locally trivial fibrations on each $\tilde{X}_{\Gamma}$
and $X_{\Gamma}$.

The following proposition gives the first properties of the strata in $\tilde{X}$.

\begin{proposition}\hspace{1in}
\begin{enumerate}
\item[(i)] If $\Gamma$ is a shifting diagram or a split diagram, then
$\tilde{X}_{\Gamma}$ is an irreducible locally closed
variety with codimension as in Tables~\ref{tab:shifting} and
\ref{tab:split}, respectively.
\item[(ii)] More generally, if $\Gamma$ is any admissible diagram with
split part $\Gamma_{\spl}$ and shifting part $\Gamma_{\shi}$, then
$\tilde{X}_{\Gamma}$ is an irreducible locally closed variety and
\[\cod(\tilde{X}_{\Gamma})=\cod(\tilde{X}_{\Gamma_{\shi}})+
\cod(\tilde{X}_{\Gamma_{\spl}}).\]
\item[(iii)] The variety $\tilde{X}_{\Gamma}$ is non-empty if and only
if the diagram $\Gamma$ is admissible.
\end{enumerate}
\end{proposition}

\begin{proof}{}
(i) follows from explicit computations using the equations of
Theorem~\ref{thm:local-equations} together with the diagrams in
Tables~\ref{tab:shifting} and \ref{tab:split}.  (ii) follows from the
observation that a general point in $\tilde{X}_{\Gamma}$ can be
described by specifying coordinates for a general ``split point'' in
$\tilde{X}_{\Gamma_{\spl}}$ together with any nonzero scaling factor
for each unmarked component of $S_{\shi}$.  For (iii), the fact that
$\tilde{X}_{\Gamma}\neq\emptyset$ implies $\Gamma$ is admissible is a
restatement of Proposition~\ref{prop:triangle-rules}.  The converse
follows from (i) and (ii).
\end{proof}

\begin{theorem}\label{thm:normalcrossings}
Let $Z\subset \tilde{X}$ be the closure of the union of all
$\tilde{X}_{\Gamma }, \Gamma \in \sS$ of codimension one.  Then $Z$ is
a divisor with normal crossings.
\end{theorem}

\begin{proof}
It suffices to prove that $Z$ is a divisor with normal crossings in
any open set $\tilde{U}$.  To do this, we first consider the
subvariety $Y\subset\tilde{U}$ defined by setting all
of the affine coordinates $u_{\alpha}=0$, $\alpha\in\cE$.  It
follows from \cite[Lemma~6.5]{BGS} that $\tilde{U}$ can be identified
with a sum of line bundles $L_1\oplus L_2\oplus L_3$ on
$Y$, and that our decomposition of $\tilde{U}$
coincides with the decomposition of this bundle constructed by
restricting the coordinate subbundles to the strata in
$Y$.  For example, the stratum of type
$AA^*CD$ ( i.e., the one corresponding to the diagram with shifting part
of type $AA^*$ and split part of type $CD$) is the restriction of the
subbundle $L_2$ (minus its zero section) to the stratum of type
$ABA^*CD$ in $Y$.

Since the coordinate subbundles of $L_1\oplus L_2\oplus L_3$ have
normal crossings in each fiber, $Z$ will be a divisor with normal
crossings provided $Z\cap Y$ is a divisor with normal
crossings in $Y$.  This can be verified by direct
computation using the techniques and rules for differentials used to
verify nonsingularity of $\tilde{X}$ \cite[Theorem~7.6]{BGS}.  There
we showed that $Y$ was a smooth variety, and singled
out coordinates whose differentials along the minimal strata formed a
basis for the cotangent space.  In fact, it is easy to show that
subsets of these differentials span bases for the cotangent spaces to
the divisors meeting along this stratum, proving that along the
minimal strata, the irreducible components of $Z\cap Y$
are smooth and have normal crossings.  Similar computations work for
the higher dimensional strata in $Y$.
\end{proof}

\begin{corollary}\label{cor:strat}
The decomposition $\{\tilde{X}_{\Gamma} \mid \Gamma \in \sS \}$ is a stratification.
\end{corollary}

\begin{proof}
Any divisor with normal crossings determines a stratification in which
the closures of the strata are precisely the multiple intersections of
the irreducible divisors.  One only needs to observe that the
decomposition of $\tilde{X}$ indexed by admissible diagrams coincides
with the stratification induced by $Z$.
\end{proof}

\section{Betti numbers}\label{s:betti}
\subsection{Zeta function and Betti numbers}
In this section we use the Hasse-Weil zeta function $Z (\tilde{X},s)$
of $\tilde{X}$ to compute the Poincar\'{e} polynomial
\[
P_{\tilde{X}}(t)=\sum_{i=0}^{24}\dim(H^i(\tilde{X} (\bbC),\bbC))\cdot t^i.
\]
In other words, we count points on $\tilde{X} (\bbF_{q})$ as a
function of $q$, and then use the relation between $Z (\tilde{X},s)$
and the ranks of the complex cohomology groups $H^{*} (\tilde{X}
(\bbC); \bbC)$ to compute the latter.  For more information about $Z$
and for the proof of Weil's conjectures, we refer to \cite{Del}.

Let $q=p^{l}$ for some prime $p$ and positive integer $l$, and let
$\bbF_{q}$ be the finite field with $q$ elements. Let $Y$ be an
$n$-dimensional nonsingular projective variety defined over
$\bbF_{q}$, and for $r\geq 1$ let $\#Y (\bbF_{q^{r}})$ be the number of
$\bbF_{q^{r}}$-rational points of $Y$.  Then the \emph{zeta function} $Z (
Y, s) $ is the formal power series defined by
\[
 Z (Y,s)= \exp (\sum_{r\geq 1} \#Y (\bbF_{q^{r}})s^{r}/r).
\]

\begin{theorem}
(cf. \cite{Del})
The zeta function is a rational function of $s$, with a factorization
\[
Z (Y,s) = \frac{P_{1} (s)\cdots P_{2n-1} (s)}{P_{0} (s)\cdots P_{2n} (s)}.
\]
If $Y$ is the reduction modulo $\bbF_{q}$ of a nonsingular
projective variety defined over $\bbC$, then the degree of $P_{i}$ is
the $i$th topological Betti number of $Y$.
\end{theorem}

The following fact is well known, but we were unable to locate a
reference; thus we provide a proof.

\begin{lemma}\label{lem:Weil}
Suppose that $\#Y (\bbF_{q^{r}})$ is given by the polynomial $\sum
_{i=0}^{n}a_{i}q^{ri}$ for all $\bbF_{q^{r}}$.  Then the rank of
$H^{k} (Y (\bbC ); \bbC) $ is $0$ if $k$ is odd, and is $a_{i}$ if $k
= 2i$.
\end{lemma}

\begin{proof}
The proof is a simple manipulation with power series:
\begin{align*}
Z (Y,s) &= \exp (\sum_{r} \#Y (\bbF_{q^{r}})s^{r}/r)\\
	&= \exp (\sum_{r}  \sum_{i} a_{i}q^{ri}s^{r}/r)\\
        &= \exp (\sum_{i} a_{i}\sum_{r} q^{ri}s^{r}/r)\\
        &= \prod_{i} \exp (a_{i}\sum_{r} (q^{i}s)^{r}/r)\\
        &= \prod_{i}\frac{1}{(1-q^{i}s)^{a_{i}}}.\\
\end{align*}
In the last step we used $-\log (1-x) = \sum x^{r}/r$.
\end{proof}

\noindent{\bf Poincar\'{e} polynomial for the flag variety.}
To illustrate how to use Lemma \ref{lem:Weil}, and to warm up for
later calculations, we compute the (well known) Poincar\'{e}
polynomial of the flag variety $\Flag$.  In order to count
points in $\Flag (\bbF_{q})$, $q=p^l$, we must count flags $V_{1}\subset
V_{2}\subset V_{3}$ in $(\bbF_{q})^{4}$.  There are
$(q^4-1)/(q-1)=1+q+q^2+q^3$ choices for $V_1$, $(q^3-1)/(q-1)=1+q+q^2$
choices for $V_2$ containing $V_1$, and $(q^2-1)/(q-1)=1+q$ choices
for $V_3$ containing $V_2$.  Thus
\[\#\Flag (\bbF_{q})= (1+q+q^2+q^3)(1+q+q^2)(1+q)=1+3q+5q^2+6q^3+5q^4+3q^5+q^6.\]
In fact, the same arguments work for any extension $\bbF_{q^{r}}$ of $\bbF_{q}$,
so by Lemma~\ref{lem:Weil},
\[P_{\Flag}(t)=1+3t^2+5t^4+6t^6+5t^8+3t^{10}+t^{12}.\]

\subsection{Counting points in $\tilde{X}$}
We fix a prime power $q=p^{l}$ for some large $p$, and
we want to count points in $\tilde{X} (\bbF_{q^{r}})$.  As above,
it will suffice to determine $\#\tilde{X} (\bbF_q)$, since our
arguments will be the same for all extensions of $\bbF_{q}$. (The
reason is that the strata are very simple varieties, usually
projective spaces).
We will count the points in the (disjoint) sets $\tilde{X}_{\Gamma} (\bbF_{q})$
separately.  A further reduction is that since
$\tilde{X}\rightarrow\Flag$ is a locally trivial fibration, the
cardinality of each $\tilde{X}_{\Gamma } (\bbF_{q})$ is the product of
$\#(\tilde{X}_{\Gamma }\cap X^{fib}) (\bbF_{q})$ with $\#\Flag
(\bbF_{q})$.  Hence we fix the flag $x_{1}\subset x_{12}\subset
x_{123}$ and count the former factor.

\begin{figure}[thb]
\begin{center}
\includegraphics[scale=0.5]{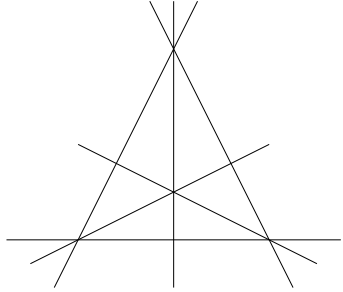}
\end{center}
\caption{The hyperplane arrangement of type ${\mathbf A}_{3}$.\label{fig:a3arrange.fig}}
\end{figure}

The counts $\#(\tilde{X}_{\Gamma }\cap X^{fib})
(\bbF_{q})$ are given in Table~\ref{tab:points}. In this table we
group the point counts by split diagrams, following the
decomposition of Proposition~\ref{prop:diagram-classification}.  For
each split diagram, we indicate the corresponding shifting diagrams
using $\bullet$ and $\circ$ as in Table~\ref{tab:shifting}.  To
construct the table we use
the following facts:
\begin{itemize}
\item For $\tilde{X}_{\Gamma}\rightarrow X_{\Gamma}$ to have a
nontrivial fiber, the diagram $\Gamma = (S,S_{\#})$ must contain
either all the edges in a hypersimplex, or at least one pair of related
triangles.  Indeed, otherwise the coordinates corresponding to every
pair of related triangles are determined, and this suffices to
determine the coodinates for all components of the graph.  This also
implies that for many strata the fiber consists of a
single point; in these cases we only need to count points in $X_{\Gamma }$.
\item According to Theorem~\ref{thm:local-equations}, any point in
$\tilde{X}_{\Gamma} $ has a neighborhood that embeds in the
product
\[
\Flag(V)\times\bbA_{\cE_{1}}\times \bbA_{\cE_{2}}\times
\bbA_{\cE_{3}}\times\bbP_{\Delta_{1}}\times \bbP_{\Delta_{2}}\times
\bbP_{\Delta_{3}}\times \prod_{\beta\in\cT} \bbP_{\beta},
\]
where $\cT$ is the set of triangular faces of the $\Delta_{i}$.
A simple computation in coordinates
proves the following:
\begin{itemize}
\item Let $\beta$ be the triangle $(a,b,c)$.  Then the subvariety of
$\bbP_{\beta}$ cut out by the linear relation $x_{a}-x_{b}+x_{c}=0$ is
isomorphic to $\bbP^{1}$, and the subvariety cut out by
$x_{a}x_{b}x_{c} = 0$ is a set of three distinct points.
\item Let $a,b,c,d,e,f$ be the edges of the tetrahedron $\Delta_{1}$.
Then the subvariety of $\bbP_{\Delta_{1}}$ cut out by the 4 linear
relations in $x_{a},\dotsc ,x_{f}$ is isomorphic to $\bbP^{2}$, and
the subvariety of this cut out by $x_{a}x_{b}x_{c}x_{d}x_{e}x_{f} = 0$
is the projective hyperplane arrangement of type ${\mathbf A}_{3}$ (Figure
\ref{fig:a3arrange.fig}).  The same is true for $\bbP_{\Delta_{3}}$.
\end{itemize}
\end{itemize}

\noindent {\bf Point counts for shifting strata.}  Most of the point
counts are straighforward.  The most difficult computations are for
the shifting strata, so we describe these in some detail.

\begin{enumerate}
\setlength{\itemsep}{.1in}
\item[] {\em Type $X_\emptyset$}  ($\circ\circ\circ$).
The map $\tilde{X}_{\Gamma }\rightarrow X_\Gamma $ is one-to-one, so
it suffices to count the points in $X_{\Gamma} = X^{\circ}$ (nondegenerate
tetrahedra).  It is easy to see that any nondegenerate tetrahedron can
be constructed by choosing a flag in general position to the fixed
flag $x_1\subset x_{12}\subset x_{123}$.  This means the number of
nondegenerate tetrahedra is the same as the number of points in the
open set of the flag variety, which is $q^{6}$.

\item[] {\em Types $A$ and $A^*$} ($\bullet\circ\circ$ \emph{and}
$\circ\circ\bullet$).  These are dual, and so have the same number of
points.  We count points for $A^*$.  We claim that the fiber of
$\tilde{X}_{\Gamma}\rightarrow X_{\Gamma}$ is a single point.  Indeed,
the coordinates of any point in $\tilde{X}_{\Gamma}$ are determined by
the projections to $\bbA_{\cE_{1}}$, $\bbA_{\cE_{2}}$, and
$\bbP_{\Delta_{3}}$, and the projection to $\bbP_{\Delta_{3}}$ is
determined (via the quadric and quartic relations from
Theorem~\ref{thm:local-equations}) by those to $\bbA_{\cE_{1}}$ and
$\bbA_{\cE_{2}}$.  Since for this stratum
$\tilde{X}_{\Gamma}\rightarrow X_{\Gamma}$ is given by projection to
$\bbA_{\cE_{1}}\times \bbA_{\cE_{2}}$, the fiber of
$\tilde{X}_{\Gamma}\rightarrow X_{\Gamma}$ is a single point.  In
fact, a similar argument shows that the projection
$\tilde{X}_{\Gamma}\rightarrow X_{\Gamma}$ is one-to-one for any
codimension-1 stratum.

To count points in $X_{\Gamma}$, note that a point of type $A^*$ is a
degenerate tetrahedron obtained by collapsing all of the planes to the
single fixed plane $x_{123}$, but is as general as possible otherwise
(see Figure~\ref{fig:shifting-div}).  Having fixed the flag
$x_1,x_{12},x_{123}$, we then have $q$ choices for $x_2$ (a point on
$x_{12}$ not equal to $x_1$), $q$ choices for $x_{13}$ (a line
containing $x_1$ and contained in $x_{123}$, but different from
$x_{12}$), $q$ choices for $x_3$ (a point on $x_{13}$ not equal to
$x_1$), $q-1$ choices for $x_{14}$ (a line containing $x_1$ and
contained in $x_{123}=x_{124}$, but different from $x_{12}$ and
$x_{13}$), and $q-1$ choices for $x_4$ (a point on $x_{14}$ different
from $x_1$, and not lying on the line $x_{23}$).  The remaining lines
are then completely determined, so there are $q^3(q-1)^2=q^5-2q^4+q^3$
points altogether.

\item[] {\em Type $B$} ($\circ\bullet\circ$). Again,
$\tilde{X}_{\Gamma}\rightarrow X_{\Gamma}$ is one-to-one
over the image of this stratum, so we count points $x\in X_{\Gamma}$.
Such an $x$ is a degenerate tetrahedron obtained by collapsing all of
the lines to the single fixed line $x_{12}$ (see
Figure~\ref{fig:shifting-div}).  There are $q$ choices for the point $x_2$,
$q-1$ choices for $x_3$, $q-2$ choices for $x_4$, $q$ choices for the
plane $x_{124}$, $q-1$ choices for $x_{134}$.  The final plane
$x_{234}$ is then determined by equation (4) of
Theorem~\ref{thm:local-equations} (this equation says
that the cross-ratio of the $4$ points on the line $x_{12}$ must be the
same as the cross-ratio of the $4$ planes containing the line $x_{12}$).
Thus, $\#\tilde{X}_{\Gamma} (\bbF_{q}) =  q^2(q-1)^2(q-2)=q^5-4q^4+5q^3-2q^2$.

\item[] {\em Types $AB$ and $BA^*$} ($\bullet\bullet\circ$ \emph{and}
$\circ\bullet\bullet$).  Again, these are dual, so it suffices to
consider $BA^*$.  This time the map $\tilde{X}_{\Gamma}\rightarrow
X_{\Gamma}$ is not one-to-one, so we must be more careful.
If we know the values of all coordinates on
$\bbA_{\cE_{1}}$ and $\bbP_{\Delta_{3}}$, then Theorem
\ref{thm:local-equations} implies that we will know all of the
coordinate values.  The points in $\bbP_{\Delta _{3}}$ from this
stratum are exactly the points in $\bbP^{2}$ not on the ${\mathbf A}_{3}$
hyperplane arrangement, and this gives $(q^{2}+q+1)-6 (q+1)+ (1\cdot 3
+ 2\cdot 4) = q^{2}-5q+6$ points.  After these values are fixed, we
can specify two points in the configuration arbitrarily, but then the
remaining two points are determined (one because it is $x_{1}$ and is
globally fixed, and one by the cross-ratio condition).  This gives a
factor of $q(q-1)$ more choices, which makes the total number of
points in this stratum ${q^{4}-{6}q^{3}+{11}q^{2}-{6}q}$.

\item[] {\em Type $AA^{*}$} ($\bullet\circ\bullet$).  To count points in
this stratum, we use a degeneration trick.  Consider the configuration
$C$ of $6$ lines in the affine plane $\bbA^{2}$ shown in
Figure~\ref{fig:collapse}.
If we linearly collapse all the points into the fixed point $x_{1}$,
then we will obtain a configuration of $6$ lines contained in the
fixed plane, and all containing the fixed point.  It is easy to show
that any point in this stratum can be obtained in this way, and that
the points obtained from two non-homothetic configurations in
$\bbA^{2}$ are distinct.

\begin{figure}[htb]
\begin{center}
\includegraphics[scale = .4]{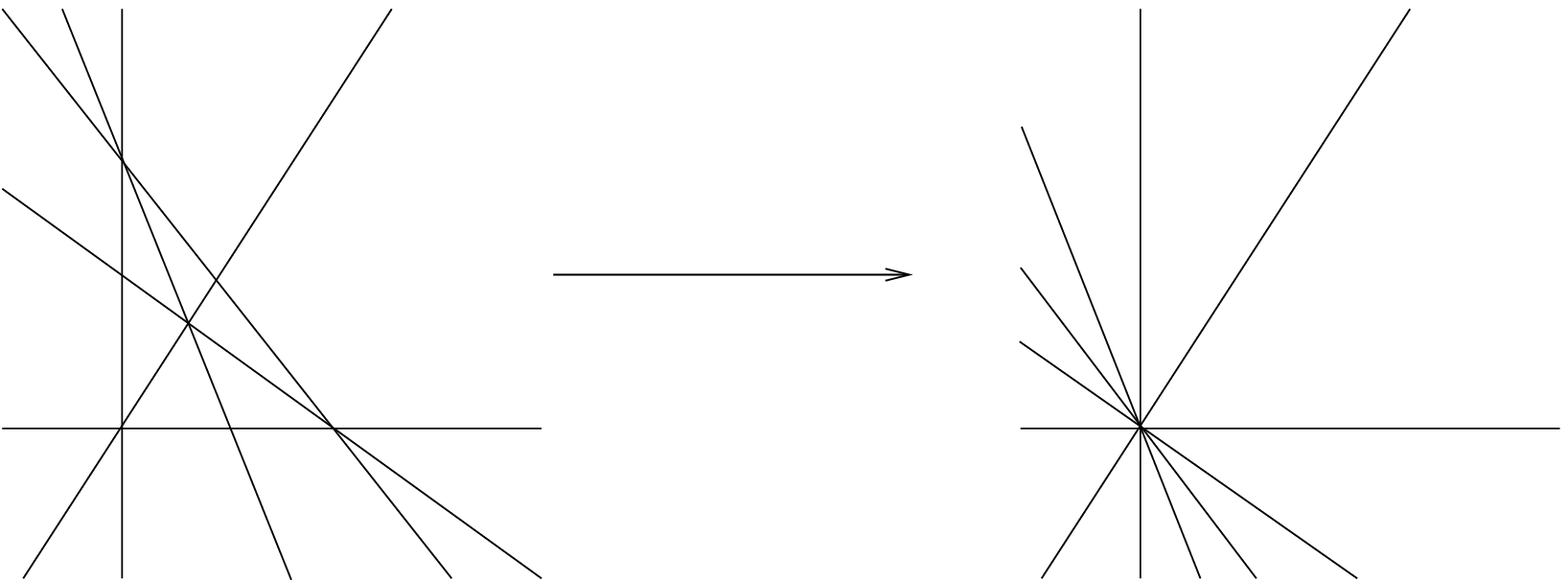}
\end{center}
\caption{\label{fig:collapse}}
\end{figure}

To build the configuration $C$, first choose the lines $x_{13}$ and
$x_{23}$ through the point $x_{1}$, and then move them slightly away
from $x_{1}$ to form a small triangle $T$.  Then $x_{4}$ can be placed
on any point in $\bbA^{2}\smallsetminus T$, and this determines the
remaining three lines.  Altogether we have $q (q-1)$ choices for
$x_{13}$ and $x_{23}$, and then $q^{2}-3q+3$ choices for $x_{4}$, so
this gives ${q^{4}-{4}q^{3}+{6}q^{2}-{3}q}$ points for this stratum.

\item[] {\em Type $ABA^{*}$} ($\bullet\bullet\bullet$). Just as for the
stratum $BA^{*}$, all values of coordinates are known if we know the
projections to $\bbP_{\Delta_{1}}$ and $\bbP_{\Delta_{3}}$.  For the
latter projection we know from the argument for type $BA^{*}$ that we
have $q^{2}-5q+6$ points in the image.  Moreover, the projection map
$\bbP_{\Delta _{1} }\times \bbP_{\Delta_{3}} \rightarrow
\bbP_{\Delta_{3}}$ induces a fibration of this stratum over
$\bbP^{2}\smallsetminus{\mathbf A}_{3}$; the quartic relations imply that the
fiber $F$ is contained in a smooth quadric $Q$ in $\bbP^{2}$.  Any
such quadric has $q+1$ points.  We want to count the points in $Q$
that miss the ${\mathbf A}_{3}$ hyperplane arrangement in $\bbP^{2}\subset
\bbP_{\Delta_{1}}$, and a computation shows that each quadric meets
the arrangement exactly in the four triple points.  Hence the number
of points in $F$ is $(q+1)-4$, which makes the total
number of points ${q^{3}-{8}q^{2}+{21}q-{18}}$.
\end{enumerate}

\noindent {\bf The remaining strata.} Counting points for the rest of
the strata is not difficult once one
determines the fibers of the maps $\tilde{X}_{\Gamma}\rightarrow
X_{\Gamma}$.  Then $\#\tilde{X}_{\Gamma} (\bbF_{q})$ is given by
multiplying $\#X_{\Gamma} (\bbF_{q})$, which can be determined by
straightforward geometric arguments, by the number of points in the
appropriate fiber.  In Table~\ref{tab:points} we indicate which strata
have nontrivial fibers, and we leave the details to the reader.  There
are three types of fibers that occur:
\[\mbox{(I)}\;\;\bbP^{2}\smallsetminus{\mathbf A}_{3},\hspace{.5in}
\mbox{(II)}\;\;\bbP^{1} \smallsetminus \{\hbox{three points} \},
\hspace{.5in}
\mbox{(III)}\;\;\bbP^{1} \smallsetminus \{\hbox{two points} \}.\]
In some cases the fiber is a product of two of these basic types; we
indicate this with an exponent.

\newcommand{\Fi}{\hfill \hbox{I}}
\newcommand{\Fii}{\hfill \hbox{II}}
\newcommand{\Fiii}{\hfill \hbox{III}}

\begin{table}[t]\caption{\label{tab:points} Point counts.}
\begin{tabular}{||c||c|l||c|l||} \hline
type & shift & points over $\bbF_q$ & shift & points over $\bbF_q$ \\
\hline\hline
$X_{\emptyset}$
&$\circ\circ\circ$   & ${q^{6}}$ 				&$\bullet\bullet\circ$ 	&  $ {q^{4}-{6}q^{3}+{11}q^{2}-{6}q} $\\
\cline{2-5}
&$\bullet\circ\circ$ & ${q^{5}-{2}q^{4}+q^{3}}$ 		&$\bullet\circ\bullet$ 	& $ {q^{4}-{4}q^{3}+{6}q^{2}-{3}q} $\\ \cline{2-5}
(1)
&$\circ\bullet\circ$ &${q^{5}-{4}q^{4}+{5}q^{3}-{2}q^{2}} $	&$\circ\bullet\bullet$ 	&  $ {q^{4}-{6}q^{3}+{11}q^{2}-{6}q} $\\ \cline{2-5}
&$\circ\circ\bullet$ & $ {q^{5}-{2}q^{4}+q^{3}} $ 		&$\bullet\bullet\bullet$ & $ {q^{3}-{8}q^{2}+{21}q-{18}} $\\
\hline\hline
$C$
&$\circ\circ\circ$   & $ {q^{5}-q^{4}} $  		&$\bullet\bullet\circ$ &  $ {q^{3}-{3}q^{2}+{2}q} $\\ \cline{2-5}
&$\bullet\circ\circ$ &$ {q^{4}-q^{3}} $  		&$\bullet\circ\bullet$ &  $ {q^{3}-{3}q^{2}+{2}q} $\\ \cline{2-5}
(4)
&$\circ\bullet\circ$ & $ {q^{4}-{3}q^{3}+{2}q^{2}} $	&$\circ\bullet\bullet$ & $ {q^{3}-{5}q^{2}+{6}q} \Fi $\\ \cline{2-5}
&$\circ\circ\bullet$ & $ {q^{4}-{3}q^{3}+{2}q^{2}} $	&$\bullet\bullet\bullet$ &  $ {q^{2}-{5}q+{6}} \Fi $\\
\hline
$C^*$
&$\circ\circ\circ$   &  $ {q^{5}-q^{4}} $		&$\bullet\bullet\circ$ &  $ {q^{3}-{5}q^{2}+{6}q}\Fi $ \\ \cline{2-5}
&$\bullet\circ\circ$ & 	$ {q^{4}-{3}q^{3}+{2}q^{2}} $	&$\bullet\circ\bullet$ & $ {q^{3}-{3}q^{2}+{2}q} $\\ \cline{2-5}
(4)
&$\circ\bullet\circ$ &	$ {q^{4}-{3}q^{3}+{2}q^{2}} $	&$\circ\bullet\bullet$ & $ {q^{3}-{3}q^{2}+{2}q} $\\ \cline{2-5}
&$\circ\circ\bullet$ & $q^4-q^3$ 			&$\bullet\bullet\bullet$ & $q^2-5q+6\Fi $\\
\hline
$D$
&$\circ\circ\circ$   & $ {q^{5}} $ 	 	&$\bullet\bullet\circ$ & $ {q^{3}-{3}q^{2}+{2}q}\Fii  $\\ \cline{2-5}
&$\bullet\circ\circ$ &$ {q^{4}-q^{3}} $ 	&$\bullet\circ\bullet$ & $ {q^{3}-{2}q^{2}+q} $ \\ \cline{2-5}
(6)
&$\circ\bullet\circ$ &$ {q^{4}-{2}q^{3}+q^{2}}$ &$\circ\bullet\bullet$ & $ {q^{3}-{3}q^{2}+{2}q}\Fii  $ \\ \cline{2-5}
&$\circ\circ\bullet$ & $ {q^{4}-q^{3}} $  	&$\bullet\bullet\bullet$ & $ {q^{2}-{4}q+{4}} \Fii^{2}$ \\
\hline
$E$
&$\circ\circ\circ$   &$ {q^{5}-{2}q^{4}+q^{3}} $ 	&$\bullet\bullet\circ$ & $ {q^{3}-{3}q^{2}+{2}q} \Fii $\\ \cline{2-5}
&$\bullet\circ\circ$ &${q^{4}-{3}q^{3}+{2}q^{2}}\Fii$   &$\bullet\circ\bullet$ & $ {q^{3}-{4}q^{2}+{4}q} \Fii^{2}$\\ \cline{2-5}
(6)
&$\circ\bullet\circ$ & $ {q^{4}-{2}q^{3}+q^{2}} $	&$\circ\bullet\bullet$ & $ {q^{3}-{3}q^{2}+{2}q} \Fii $\\ \cline{2-5}
&$\circ\circ\bullet$ &$ {q^{4}-{3}q^{3}+{2}q^{2}}\Fii$ 	&$\bullet\bullet\bullet$ & $ {q^{2}-{4}q+{4}} \Fii^{2}$\\
\hline\hline
$CC_{nop}^*$
&$\circ\circ\circ$   &$ {q^{4}-q^{3}} $  &$\bullet\bullet\circ$ & $ {q^{2}-{2}q}\Fii  $\\ \cline{2-5}
&$\bullet\circ\circ$ &$ {q^{3}-q^{2}} $  &$\bullet\circ\bullet$ & $ {q^{2}-q} $\\ \cline{2-5}
(12)
&$\circ\bullet\circ$ &${q^{3}-{2}q^{2}}\Fii $  &$\circ\bullet\bullet$ & $ {q^{2}-{2}q} \Fii $ \\ \cline{2-5}
&$\circ\circ\bullet$ &$ {q^{3}-q^{2}} $        &$\bullet\bullet\bullet$ & $ {q-{2}} \Fii $\\
\hline
$CC_{op}^*$
&$\circ\circ\circ$   & $q^4-2q^3\Fii$		&$\bullet\bullet\circ$ & $q^2-2q\Fii$\\ \cline{2-5}
&$\bullet\circ\circ$ & $q^3-2q^2\Fii$		&$\bullet\circ\bullet$ & $q^2-2q\Fii$\\ \cline{2-5}
(4)
&$\circ\bullet\circ$ & $q^3-2q^2\Fii$ 		&$\circ\bullet\bullet$ & $q^2-2q\Fii$\\ \cline{2-5}
&$\circ\circ\bullet$ & $q^3-2q^2\Fii$ 		&$\bullet\bullet\bullet$ & $q-2\Fii$\\
\hline
$CD$
&$\circ\circ\circ$   &$ {q^{4}} $ 	&$\bullet\bullet\circ$ & $ {q^{2}-q} $\\ \cline{2-5}
&$\bullet\circ\circ$ &$ {q^{3}} $  	&$\bullet\circ\bullet$ & $ {q^{2}-q} $\\ \cline{2-5}
(12)
&$\circ\bullet\circ$ &$ {q^{3}-q^{2}} $ &$\circ\bullet\bullet$ & $ {q^{2}-{2}q} \Fii $\\ \cline{2-5}
&$\circ\circ\bullet$ &$ {q^{3}-q^{2}} $ &$\bullet\bullet\bullet$ & $ {q-{2}} \Fii $\\
\hline
$C^*D$
&$\circ\circ\circ$   & $q^4$		&$\bullet\bullet\circ$ & $q^2-2q\Fii $ \\ \cline{2-5}
&$\bullet\circ\circ$ & $q^3-q^2$	&$\bullet\circ\bullet$ & $q^2-q$ \\ \cline{2-5}
(12)
&$\circ\bullet\circ$ & $q^3-q^2$	&$\circ\bullet\bullet$ & $q^2-q$ \\ \cline{2-5}
&$\circ\circ\bullet$ & $q^3$		&$\bullet\bullet\bullet$ & $q-2\Fii$\\
\hline
$CE$
&$\circ\circ\circ$   &$ {q^{4}-q^{3}} $ 	&$\bullet\bullet\circ$ &$ {q^{2}-q} $ \\ \cline{2-5}
&$\bullet\circ\circ$ &$ {q^{3}-q^{2}} $ 	&$\bullet\circ\bullet$ & $ {q^{2}-{2}q}\Fii  $\\ \cline{2-5}
(12)
&$\circ\bullet\circ$ &$ {q^{3}-q^{2}} $ 	&$\circ\bullet\bullet$ & $ {q^{2}-{2}q} \Fii $\\ \cline{2-5}
&$\circ\circ\bullet$ &$ {q^{3}-{2}q^{2}}\Fii$ 	&$\bullet\bullet\bullet$ & $ {q-{2}} \Fii $\\
\hline
\end{tabular}
\end{table}

\begin{table}[ht]
\begin{tabular}{||c||c|l||c|l||} \hline
type & shift & points over $\bbF_q$  & shift & points over $\bbF_q$ \\
\hline
$C^*E$
&$\circ\circ\circ$   & $q^4-q^3$	&$\bullet\bullet\circ$ & $q^2-2q\Fii $\\ \cline{2-5}
&$\bullet\circ\circ$ & $q^3-2q^2\Fii $	&$\bullet\circ\bullet$ & $q^2-2q\Fii $\\ \cline{2-5}
(12)
&$\circ\bullet\circ$ & $q^3-q^2$	&$\circ\bullet\bullet$ & $q^2-q$\\ \cline{2-5}
&$\circ\circ\bullet$ & $q^3-q^2$	&$\bullet\bullet\bullet$ & $q-2\Fii $\\
\hline
$DD_{op}$
&$\circ\circ\circ$   & $ {q^{4}} $	&$\bullet\bullet\circ$ &$ {q^{2}-q}\Fiii  $ \\ \cline{2-5}
&$\bullet\circ\circ$ & $ {q^{3}} $ 	&$\bullet\circ\bullet$ & $ {q^{2}} $\\ \cline{2-5}
(3)
&$\circ\bullet\circ$ &$q^3-q^2\Fiii$    &$\circ\bullet\bullet$ & $ {q^{2}-q} \Fiii $\\ \cline{2-5}
&$\circ\circ\bullet$ &$ {q^{3}} $ 	&$\bullet\bullet\bullet$ & $ {q-1} \Fiii $\\
\hline
$DE$
&$\circ\circ\circ$   &$ {q^{4}-q^{3}} $  	&$\bullet\bullet\circ$ & $ {q^{2}-{2}q}\Fii  $ \\ \cline{2-5}
&$\bullet\circ\circ$ &$ {q^{3}-{2}q^{2}}\Fii$  	&$\bullet\circ\bullet$ & $ {q^{2}-{2}q}\Fii  $\\ \cline{2-5}
(6)
&$\circ\bullet\circ$ & $ {q^{3}-q^{2}} $ 	&$\circ\bullet\bullet$ & $ {q^{2}-q} $\\ \cline{2-5}
&$\circ\circ\bullet$ &$ {q^{3}-q^{2}} $ 	&$\bullet\bullet\bullet$ & $ {q-{2}}\Fii  $\\
\hline
$DE^*$
&$\circ\circ\circ$   & $q^4-q^3$	&$\bullet\bullet\circ$ & $q^2-q$\\ \cline{2-5}
&$\bullet\circ\circ$ & $q^3-q^2$	&$\bullet\circ\bullet$ & $q^2-2q\Fii $\\ \cline{2-5}
(6)
&$\circ\bullet\circ$ & $q^3-q^2$	&$\circ\bullet\bullet$ & $q^2-2q\Fii $\\ \cline{2-5}
&$\circ\circ\bullet$ & $q^3-2q^2\Fii $	&$\bullet\bullet\bullet$ & $q-2\Fii $\\

\hline\hline
$CC_{op}^*D$
&$\circ\circ\circ$   &$ {q^{3}} $ 	&$\bullet\bullet\circ$ &$ {q} $ \\ \cline{2-5}
&$\bullet\circ\circ$ & $ {q^{2}} $	&$\bullet\circ\bullet$ &$ {q} $ \\ \cline{2-5}
(12)
&$\circ\bullet\circ$ &$ {q^{2}} $ 	&$\circ\bullet\bullet$ &$ {q} $ \\ \cline{2-5}
&$\circ\circ\bullet$ & $ {q^{2}} $	&$\bullet\bullet\bullet$ & $ 1 $\\
\hline
$CC_{nop}^*D$
&$\circ\circ\circ$   & $ {q^{3}} $	&$\bullet\bullet\circ$ & $ {q} $\\ \cline{2-5}
&$\bullet\circ\circ$ & $ {q^{2}} $	&$\bullet\circ\bullet$ & $ {q} $\\ \cline{2-5}
(12)
&$\circ\bullet\circ$ & $ {q^{2}} $	&$\circ\bullet\bullet$ & $ {q} $\\ \cline{2-5}
&$\circ\circ\bullet$ & $ {q^{2}} $	&$\bullet\bullet\bullet$ & $ 1 $\\
\hline
$CC^*E$
&$\circ\circ\circ$   & $ {q^{3}} $	&$\bullet\bullet\circ$ & $ {q} $ \\ \cline{2-5}
&$\bullet\circ\circ$ & $ {q^{2}} $ 	&$\bullet\circ\bullet$ & $ {q} $  \\ \cline{2-5}
(24)
&$\circ\bullet\circ$ &$ {q^{2}} $ 	&$\circ\bullet\bullet$ & $ {q} $  \\ \cline{2-5}
&$\circ\circ\bullet$ & $ {q^{2}} $	&$\bullet\bullet\bullet$ & $ 1 $\\
\hline
$CDE$
&$\circ\circ\circ$   &$ {q^{3}} $ &$\bullet\bullet\circ$ & $ {q} $\\ \cline{2-5}
&$\bullet\circ\circ$ &$ {q^{2}} $ &$\bullet\circ\bullet$ & $ {q} $\\ \cline{2-5}
(12)
&$\circ\bullet\circ$ &$ {q^{2}} $ &$\circ\bullet\bullet$ & $ {q} $\\ \cline{2-5}
&$\circ\circ\bullet$ &$ {q^{2}} $ &$\bullet\bullet\bullet$ & $ 1 $\\
\hline
$C^*DE$
&$\circ\circ\circ$   & $q^3$	&$\bullet\bullet\circ$ & $q$\\ \cline{2-5}
&$\bullet\circ\circ$ & $q^2$	&$\bullet\circ\bullet$ & $q$\\ \cline{2-5}
(12)
&$\circ\bullet\circ$ & $q^2$	&$\circ\bullet\bullet$ & $q$\\ \cline{2-5}
&$\circ\circ\bullet$ & $q^2$	&$\bullet\bullet\bullet$ & $1$\\
\hline
$DD_{op}E$
&$\circ\circ\circ$   &$ {q^{3}} $ 	&$\bullet\bullet\circ$ & $ {q} $\\ \cline{2-5}
&$\bullet\circ\circ$ &$ {q^{2}} $  	&$\bullet\circ\bullet$ &$ {q} $ \\ \cline{2-5}
(6)
&$\circ\bullet\circ$ &$ {q^{2}} $ 	&$\circ\bullet\bullet$ & $ {q} $\\ \cline{2-5}
&$\circ\circ\bullet$ &$ {q^{2}} $  	&$\bullet\bullet\bullet$ & $ 1 $\\
\hline\hline
\end{tabular}
\end{table}
\clearpage

\begin{theorem}\label{thm:ppoly}
The Poincar\'{e} polynomial for $\tilde{X}$ is
$P_{\tilde{X}}(t)=P_{fib}(t)\cdot P_{\Flag}(t)$ where
\[P_{fib}(t)=1+23t^2+114t^4+189t^6+114t^8+23t^{10}+t^{12}\]
and
\[P_{\Flag}(t)=1+3t^2+5t^4+6t^6+5t^8+3t^{10}+t^{12}.\]
\end{theorem}

\begin{proof}
Adding up the total numbers of points in each stratum (of the fiber)
gives
\begin{equation}\label{eq:count}
1+23q+114q^2+189q^3+114q^4+23q^5+q^6.
\end{equation}
The result then follows from Lemma \ref{lem:Weil} and the remark that
counting points over $\bbF_{q^{r}}$ for $r>1$ yields \eqref{eq:count}
with $q$ replaced with $q^{r}$.
\end{proof}

\section{Cohomology}\label{s:ring}
In this section we compute the rational cohomology ring
$H^*(\tilde{X};\bbQ)$.  By Theorem~\ref{thm:ppoly}, this ring is
trivial in odd degrees.  In degree $2$, there are $23$ classes
obtained by taking the Poincar\'{e} duals of the irreducible divisors
coming from our stratification.  These correspond to the diagrams of
type $A$, $B$, $A^*$, $C_i$, $C_i^*$, $D_{ij}$, and $E_{ij}$ (and the
degenerate tetrahedra shown in Figures~\ref{fig:shifting-div}
and~\ref{fig:split-div}).  There are additional degree $2$ classes
obtained by taking Poincar\'{e} duals of the special position
divisors.  Recall that each of these $14$ divisors corresponds
to a proper nonempty subset $I\subset\sets{4}$, and is obtained by
fixing a codimension $\card{I}$ subspace $V\subset\bbC^4$ and
requiring the $I$th face of a tetrahedron to have nonempty
intersection with $\bbP(V)$ in $\bbP^3$.  We denote this divisor by
$Y_{I} = Y_I(V)$, and note that its class is independent of the
choice of $V$.  We now prove that the cohomology ring has the
following description:

\begin{theorem}\label{thm:ringpres}
The cohomology ring $H^*(\tilde{X};\bbQ)$ is generated in degree
$2$ by the Poincar\'e duals of the the $23$ divisors from
Figures~\ref{fig:shifting-div} and~\ref{fig:split-div}, and the
special position divisors $Y_{I}$.  If we denote these
dual classes by $a$, $b$, $a^*$, $c_i$, $c_i^*$, $d_{ij}$, $e_{ij}$,
and $y_i$, $y_{ij}$, $y_{ijk}$ (where the indices are regarded as {\em
unordered} subsets of $\sets{4}$), then the ideal of relations is
generated by the following polynomials:
\begin{enumerate}
\item[(i)] $y_{ij}-y_i-y_j+a+c_k+c_l+d_{ij}+e_{ij}$,\\
$y_{ijk}-y_{ij}-y_{ik}+y_i+b+c_i+c_l^*+d_{jk}+e_{jk}+e_{il}+e_{jl}+e_{kl}$,\\
$y_{ij}-y_{ijk}-y_{ijl}+a^*+c_i^*+c_j^*+d_{kl}+e_{ij}$.\\
\item[(ii)]  $c_ic_j,\;\; c_i^*c_j^*,\;\; c_id_{ij},\;\;c_i^*d_{ij},\;\;
c_ie_{ij},\;\; c_i^*e_{jk},\;\; d_{ij}e_{ik},\;\; e_{ij}e_{ik},\;\;
e_{ij}e_{kl}$.\\
\item[(iii)] $a(y_i-y_j),\; b(y_{ij}-y_{ik}),\; a^*(y_{ijk}-y_{ijl}),\\
c_i(y_j-y_k),\; c_i(y_{ij}-y_{ik}),\; c_i^*(y_{jk}-y_{jl}),\;
c_i^*(y_{ijk}-y_{ijl}),\\
d_{ij}(y_i-y_j),\; d_{ij}(y_{ik}-y_{jk}),\; d_{ij}(y_{ikl}-y_{jkl}),\\
e_{ij}(y_i-y_j),\; e_{ij}(y_{kl}-y_{ik}),\;
e_{ij}(y_{ijk}-y_{ijl})$.\\
\item[(iv)] $y_i^2+y_{ij}^2+y_{ijk}^2-y_iy_{ij}-y_{ij}y_{ijk}$,\\
$y_{ij}^3-2y_iy_{ij}^2+2y_i^2y_{ij}$,\\
$y_i^4$.
\end{enumerate}
\end{theorem}

\subsection{Relations}\label{ss:relations}

To prove Theorem~\ref{thm:ringpres}, we first need to show that the
relations (i)-(iv) all hold in $H^*(\tilde{X};\bbQ)$.

The linear relations (i) follow from rational equivalences
corresponding to certain cross-ratios.  Let $\tilde{x}\in\tilde{X}$,
and recall that the images of $\tilde{x}$ in $X$ and $\Grass_I$ are
denoted by $x$ and $x_I$, respectively.  Let $V,V'\subset\bbC^4$ be
hyperplanes such that $V$, $V'$, and $x_{ij}$ are all in general
position (in $\bbP^3$).  Since the four points $\alpha=x_i$,
$\beta=x_j$, $\gamma=x_{ij}\cap V$, $\delta=x_{ij}\cap V'$ in $\bbP^3$
all lie on the line $x_{ij}$, the cross-ratio
\[\frac{\alpha\beta:\gamma\delta}{\alpha\gamma:\beta\delta}\]
defines a rational function on $\tilde{X}$.  The numerator
vanishes on the divisors $A$, $C_k$, $C_l$, $D_{ij}$, $E_{ij}$, and
$Y_{ij}(V\cap V')$.  The denominator vanishes on $Y_i(V)$ and
$Y_j(V')$.  A calculation in local coordinates shows that the order of vanishing on
these divisors is always $1$, giving the relation
\[y_{ij}-y_i-y_j+a+c_k+c_l+d_{ij}+e_{ij}.\]
For the second relation in (i), let $V,V'\subset\bbC^4$ be
$2$-dimensional subspaces such that $V\cap V'$ is $1$-dimensional
(equivalently, $V+V'$ is $3$-dimensional).  For suitably general $V$
and $V'$ one obtains four points $\alpha=x_{ij}\cap(V+V')$,
$\beta=x_{ik}\cap(V+V')$, $\gamma=x_{ijk}\cap V$, $\delta=x_{ijk}\cap
V'$ in $\bbP^3$, which all lie on the line $x_{ijk}\cap(V+V')$.  The
resulting cross-ratio defines a rational function on $\tilde{X}$ that
gives the second linear relation.  The third relation in (i) is the
dual of the first.

The monomial relations (ii) follow from the fact that the
corresponding pairs of divisors are all disjoint in $\tilde{X}$.  This
can be verified using the rules of
Proposition~\ref{prop:triangle-rules} and the diagrams in
Table~\ref{tab:split}.  For example, if a point $\tilde{x}$ were in
both of the divisors $C_i$ and $C_j$, then all of the edges in the
first component of the $S_{\#}(\tilde{x})$ diagram would have to be
bold (using the second rule of \ref{prop:triangle-rules}), which would
violate the first rule of \ref{prop:triangle-rules}.  Thus the product
$c_ic_j$ must be zero.  The remaining products are similar.

The relations (iii) follow from the observation that if the $I$th and
$J$th planes of a complete tetrahedron $\tilde{x}$ coincide along a
split divisor or a shifting divisor, and one of these planes is in
special position, then so is the other.  For example, if $\tilde{x}$
is in the divisor $A$ and also in $Y_i(V)$, then all of the points of
$x$ coincide, and the $i$th point $x_i$ lies on the hyperplane
$\bbP(V)\subset\bbP^3$.  It follows that the $j$th point $x_j$ must
also lie on this hyperplane, hence $\tilde{x}\in Y_j(V)$.  Since the
divisor $A$ meets both $Y_i(V)$ and $Y_j(V)$ transversally, we have
$ay_i=ay_j$ in $H^*(\tilde{X};\bbQ)$.  The other relations in (iii)
are similar.

The relations (iv) are induced from relations in the flag variety
$\Flag$.  The fibration $f:\tilde{X}\rightarrow\Flag$ induces a
homomorphism
\[H^{*}(\Flag;\bbQ)\stackrel{f^*}{\rightarrow}
H^{*}(\tilde{X};\bbQ),\]
and the classes $y_1,y_{12},y_{123}$ are the
images of the (Poincar\'{e} duals of the) usual complex codimension-one
Schubert cycles in $H_{10}(\Flag;\bbQ)$.  The relations among these
Schubert cycles are well-known (see, e.g., \cite{BGG} or \cite{D}),
and the relations in (iv) are the induced relations in
$H^{*}(\tilde{X};\bbQ)$, together with their images under the
symmetric group action.  (In fact, one gets the same ideal using only
the three relations involving just $y_1$, $y_{12}$, $y_{123}$.)

\subsection{Sufficiency of the relations}

Let $R^*$ denote the graded quotient ring
\[\bbZ[a,a^*,b,c_i,c_i^*,d_{ij},e_{ij},y_i,y_{ij},y_{ijk}]/I\]
where $I$ is the ideal generated by the relations (i)-(iv) of
Theorem~\ref{thm:ringpres}.  Since all of these relations hold in
$H^*(\tilde{X};\bbQ)$, there exists a graded ring homomorphism
$\phi:R^*\otimes\bbQ\rightarrow H^{2*}(\tilde{X};\bbQ)$.  To
complete the proof of Theorem~\ref{thm:ringpres}, we need to show
that $\phi$ is an isomorphism.  Using the software package Macaulay2
\cite{macaulay2}, it can be shown that:
\begin{enumerate}
\item[(a)] {\em The Hilbert series for $R^*\otimes\bbQ$ (and for
$R^*\otimes\bbF_2$) is the same as the Poincar\'{e} polynomial for
$\tilde{X}$ given in Theorem~\ref{thm:ppoly}.  In other words, 
\[\dim_{\bbQ}
R^i\otimes\bbQ=\dim_{\bbF_2}R^i\otimes\bbF_2= \dim_{\bbQ}
H^{2i}(\tilde{X};\bbQ)\]
for all $i$.}  
The Hilbert series for $R^*$ can be computed both over $\bbF_2$ and
over $\bbQ$ using Macaulay, and is the same in both cases.  The former
calculation can be done directly with the given presentation (removing
the obvious redundancies among the relations); but over $\bbQ$, in
order to speed up the computation we had to first reduce the number
of variables by using the linear relations to eliminate the classes
$y_i$ for $i\neq 1$, $y_{ij}$ for $ij\neq 12$, and $y_{ijk}$ for
$ijk\neq 123$. 
\item[(b)] {\em The multiplication pairing
\[R^i\otimes R^{12-i}\rightarrow R^{12}\cong\bbZ\oplus\{\mbox{\em odd
torsion}\}\]
is nondegenerate over $\bbF_2$ and $\bbQ$ for $0\leq i\leq 12$.}  We used
Macaulay to show this over $\bbF_2$.  The result over $\bbQ$ then
follows from the Hilbert series calculation. 
\item[(c)] {\em The element $\phi(y_1^3y_{12}^2y_{123}aba^*c_1c_1^*d_{23})$
is nonzero in $H^{24}(\tilde{X};\bbQ)$.} The intersection of the
special position divisors correponds to fixing the flag
$(x_1,x_{12},x_{123})$, and the intersection of the remaining divisors
determines a point in the fiber of $\tilde{X}\rightarrow\Flag$ over
this fixed flag.  Thus, the indicated element is the (dual of) the
class of a point in $\tilde{X}$.
\end{enumerate}
The homomorphism $\phi$ is injective by (b) and (c) and hence, an isomorphism by (a).

\bibliographystyle{amsplain}
\bibliography{betti-ring}
\end{document}